\documentclass{tac}

\usepackage{amsmath,amsfonts,amssymb}
\usepackage{setspace}
\usepackage[all,cmtip,2cell]{xy}
\newdir{ >}{{}*!/-5pt/\dir{>}}
\UseTwocells
\usepackage{tikz}
\usepackage{ifthen}
\usepackage{graphicx}
\usepackage{MnSymbol}
\usepackage{adjustbox}

\usepackage[colorlinks=true]{hyperref}
\hypersetup{allcolors=[rgb]{0.1,0.1,0.4}}

\author{Florian De Leger}

\thanks{The author is supported by RVO:67985840 and Preamium Academiae of M. Markl}

\address{Mathematical Institute of the Academy \\ \v{Z}itn\'a 25, 115~67 Prague 1, Czech Republic}

\title[Cofinal morphism of polynomial monads and double delooping]{Cofinal morphism of polynomial monads and double delooping}

\copyrightyear{2023}

\keywords{operads, polynomial monads, delooping}
\amsclass{18D20, 18D50, 55P48}

\eaddress{de-leger@math.cas.cz}

\newtheorem{thm}{Theorem}
\newtheorem{lem}{Lemma}

\let\pf\proof
\let\pfof\proofof
\let\epf\endproof
\mathrmdef{Alg}
\mathrmdef{Cat}
\mathrmdef{CAT}
\mathrmdef{colim}
\mathrmdef{Int}
\mathrmdef{Set}
\mathrmdef{SSet}
\mathbfdef{Poly}

\newcommand{\F}{\mbox{$\mathcal F$}}
\newcommand{\btimes}{\mathbin{\rotatebox{45}{$\boxtimes$}}}
\newcommand{\bempty}{\mathbin{\rotatebox{45}{$\square$}}}
\newcommand{\nfourdots}{\mathbin{\hspace{.8pt}\raisebox{-2pt}{$\diamonddots$}}}
\newcommand{\nfivedots}{\mathbin{\hspace{.8pt}\raisebox{-2pt}{$\fivedots$}}}
\newcommand{\nspecialthreedots}{\mathbin{\hspace{.8pt}\raisebox{-2pt}{\adjustbox{trim=0pt 1pt 0pt -1pt,clip,raise=1pt}{$\diamonddots$}}}}
\newcommand{\nspecialfourdots}{\mathbin{\hspace{.8pt}\raisebox{-2pt}{\adjustbox{trim=0pt 1pt 0pt -1pt,clip,raise=1pt}{$\fivedots$}}}}
\newcommand{\wbimodbfive}{\mathrm{IBimod}_{\nfivedots}}
\newcommand{\nopbfive}{\mathrm{NOp}_{\btimes}}

\begin{document}
	
\maketitle

\begin{abstract} Using the theory of internal algebras classifiers developed by Batanin and Berger, we construct a morphism of polynomial monads which we prove is homotopically cofinal. We then describe how this result constitutes the main conceptual argument for a categorical direct double delooping proof of the Turchin-Dwyer-Hess theorem concerning the explicit double delooping of spaces of long knots.
\end{abstract}

\section{Introduction}

In \cite{batanindeleger}, Batanin and the author extended some classical results of the homotopy theory of small categories and presheaves to polynomial monads and their algebras. A notion of \emph{homopically cofinal} morphism between polynomial monad was given, which extends the notion of \emph{homotopy left cofinal} functors \cite{hirschhorn}. In this present work, we will construct a homotopically cofinal morphism of polynomial monads and use it to give a direct double delooping proof of the Turchin-Dwyer-Hess theorem \cite{dwyerhess,turchin}. Recall that a non-symmetric operad $\mathcal{O}$ is \emph{multiplicative} if it is equipped with an operadic map from the terminal non-symmetric operad $\mathcal{A}ss$ to $\mathcal{O}$. Such a map endows the collection $(\mathcal{O}_n)_{n \geq 0}$ with a structure of cosimplicial object \cite{turchin}, which we will write $\mathcal{O}^\bullet$. The Turchin-Dwyer-Hess theorem can be stated as follows:
\thm\label{theoremturchindwyerhess}
If $\mathcal{O}$ is a simplicial multiplicative operad such that $\mathcal{O}_0$ and $\mathcal{O}_1$ are contractible, then there is a weak equivalence between simplicial sets
\[
\Omega^2 Map_{NOp}(\mathcal{A}ss,u^*(\mathcal{O})) \sim Tot(\mathcal{O}^\bullet),
\]
where $\Omega$ is the loop space functor, $Map_{NOp}(-,-)$ is the homotopy mapping space in the category of simplicial non-symmetric operads, $u^*$ is the forgetful functor from multiplicative to non-symmetric operads and $Tot(-)$ is the homotopy totalization.
\endthm
This result is remarkable especially because of an earlier result from Sinha \cite{sinha1} which connects the space of \emph{long knots modulo immersions} \cite{dwyerhess} with the totalization of the Kontsevich operad.
Due to its importance in algebraic geometry, many proofs of the Turchin-Dwyer-Hess theorem, and generalisation to higher dimensions, exist in the literature \cite{batanindeleger,boavida,ducoulombier,ducoulombierturchin,dwyerhess,turchin}. We will use here in particular the theory of polynomial monads as we did in \cite{batanindeleger}. Some of the results of \cite{batanindeleger} are unfortunately difficult to apply in order to prove higher delooping results. For example, in the delooping theorem \cite[Theorem 8.1]{batanindeleger}, a technical condition of \emph{homotopically cofinal square} needs to be satisfied. Such condition is sometimes difficult to check and involves complicated combinatorics. That is why, in this paper, we present new results and suggest new techniques which we believe are easier to apply in general. In fact, we have recently applied in \cite{delegergrego} the techniques presented here in order to extend the results from Turchin and Dwyer-Hess to iterations of the Baez-Dolan construction \cite{baezdolan}.


This paper is structured as follows. We will start by recalling the notions of polynomial monads, internal algebra classifiers and homotopically cofinal maps in Section \ref{chaptercofinality}. Our presentation is strongly inspired from \cite{batanindeleger}. In Section \ref{sectioncofinality}, we will construct a polynomial monad $\wbimodbfive$ whose algebras are quintuples $(A,B,C,D,E)$ where $A$ and $B$ are non-symmetric operads, $C$ and $D$ are pointed $A-B$-bimodules and $E$ is a pointed infinitesimal $C-D$-bimodule. This data can be organised in the following diagram:
\[
\xymatrix{
	& E \ar@{.}[rr] && **[r] \text{pointed infinitesimal $C-D$-bimodule} \\
	C \ar@{.}[rr] && D \ar@{.}[r] & **[r] \text{pointed $A-B$-bimodules} \\
	A \ar@{.}[rr] && B \ar@{.}[r] & **[r] \text{non-symmetric operads}
}
\]
The notions of pointed bimodules and pointed infinitesimal bimodules are introduced in Subsection \ref{subsectionpointedinfinitesimalbimodule}. The reason why $E$ is called \emph{infinitesimal} is that it coincides with the classical notion of infinitesimal bimodules (called \emph{weak bimodules} in \cite{turchin}) when $A=B=C=D$. It is not hard to see how this construction can be generalised to higher dimensions: we would have a tower of pairs of elements as in the previous diagram, where each pair are some kind of bimodules over the pair below, and with a single element sitting on top. We will also construct a polynomial $\nopbfive$
whose algebras are commutative diagrams
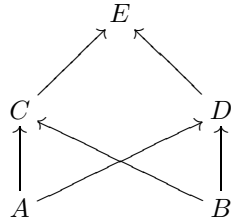
\[
\xymatrix{
	& E \\
	C \ar[ru] && D \ar[lu] \\
	A \ar[u] \ar[rru] && B \ar[u] \ar[llu]
}
\]
of non-symmetric operads. The non-symmetric operad $E$ should be seen as sitting on top of a ``circle'' of non-symmetric operads. The objective of Section \ref{sectioncofinality} is to prove Theorem \ref{theoremcofinality} which states that there is a homotopically cofinal morphism of polynomial monads from $\wbimodbfive$ to $\nopbfive$. In Section \ref{chapterapplication}, we will show how this cofinality result can be applied to the Turchin-Dwyer-Hess theorem. We will construct circles $S^1$ in any simplicial model category $M$ and show that, if $M$ is left proper, computing the double loop space of mapping spaces in $M$ reduces to studying mapping spaces in the comma category $S^1 / M$. Because $S^1$ is not a set theoretical object but is homotopical by nature, we can not directly apply our cofinality result. However, we can look at the categories of algebras of our previously defined polynomial monads as fibred categories, and prove that our cofinality result can be applied to the fibres. The last technical step is about the ``pointedness'' of infinitesimal bimodules. We will prove that mapping spaces in the categories of infinitesimal bimodule and pointed infinitesimal bimodules are weakly equivalent.

\section{Homotopy theory for polynomial monads}\label{chaptercofinality}

\subsection{Polynomial monads and their algebras}\label{sectionpolymon}

\begin{definition}[\cite{bataninberger,gambinokock}] A \emph{(finitary) polynomial} $P$ is a diagram in $\Set$, the category of sets, of the form 
	\[
	\xymatrix{
		J & E \ar[l]_-s \ar[r]^-p & B \ar[r]^-t & I
	}
	\]
	where $p^{-1}(b)$ is finite for any $b\in B$. \end{definition} 
Each polynomial $P$ generates a functor called {\it polynomial functor} between functor categories $$\underline{P}:\Set^J \to \Set^I$$ which is defined as the composite functor
\[
\xymatrix{
	\Set^J \ar[r]^-{s^*} & \Set^E  \ar[r]^-{p_*} & \Set^B \ar[r]^-{t_!} & \Set^I
}
\]
where we consider the sets $J,E,B, I$ as discrete categories, $s^*$ is the restriction functor and $p_*$ and $t_!$ are the right and left Kan extensions correspondingly.
Explicitly the functor $\underline{P}$ is given by the formula
\begin{equation}\label{PPP}
	\underline{P}(X)_i = \coprod_{b\in t^{-1}(i)} \prod_{e\in p^{-1}(b)} X_{s(e)}.
\end{equation}

A \emph{cartesian morphism} between polynomial functors from $\Set^J$ to $\Set^I$ is a natural transformation between the functors such that each naturality square is a pullback. 
Composition of finitary polynomial functors is again a finitary polynomial functor. Sets, finitary polynomial functors and their cartesian morphisms form a $2$-category  $\Poly$.
\begin{definition}
	A \emph{(finitary) polynomial monad} is a monad in the $2$-category $\Poly$.
\end{definition}

Polynomial monads form a category. Explicitly, a morphism of polynomial monads $f: S \to T$ can be described \cite{bataninberger} as a commutative diagram
\begin{equation}\label{morphismpolymon}
\xymatrix{
	J \ar[d]_\phi & D \ar[l]_v \ar[r]^q \ar[d]_\pi & C \ar[r]^u \ar[d]^\psi & J \ar[d]^\phi \\
	I & E \ar[l]^s \ar[r]_p & B \ar[r]_t & I
}
\end{equation}
in which the middle square is a pullback, where the top horizontal line is the polynomial generating $S$ and the bottom horizontal line is the polynomial generating $T$.

For a polynomial monad $T$ given by
\begin{equation*}
\xymatrix{
	I & E \ar[l]_-s \ar[r]^-p & B \ar[r]^-t & I
}
\end{equation*}
we will call $I$ the set of \emph{colours} of $T$, $B$ the set of \emph{operations}, $E$ the set of \emph{marked operations}, the map $t$ \emph{target} and the map $s$ \emph{source}. The map $p$ will be called the \emph{middle map} of $T$.

Examples of polynomial monads are small categories, the free monoid monad and the free monad $\mathrm{NOp}$ for non-symmetric operads \cite{batanindeleger}. Let us describe $\mathrm{NOp}$ again, as it is an important example for us. Recall that a \emph{non-symmetric operad} $A$ in a symmetric monoidal category $(\mathcal{E},\otimes,I)$ is given by
\begin{itemize}
	\item an object $A_n$ in $\mathcal{E}$ for all integers $n \geq 0$
	\item a morphism $\epsilon: I \to A_1$ called \emph{unit}
	\item for any $n \geq 0$, $m_1,\ldots,m_n \geq 0$ with $m := m_1 + \ldots m_n$, a map
	\begin{equation}\label{mapnonsymmetricoperad}
	\mu : A_n \otimes \bigotimes_{j=1}^n A_{m_j} \to A_m
	\end{equation}
	called \emph{multiplication},
\end{itemize}
such that the usual associativity and unity conditions are satisfied.

\begin{example} \label{example3}
	The polynomial for $\mathrm{NOp}$ is given by
	\[
	\xymatrix{
		\mathbb{N} & PTr^{*} \ar[l]_-s \ar[r]^-{p} & PTr \ar[r]^-t & \mathbb{N}
	}
	\] 
	where $PTr$ and $PTr^*$ are the sets of isomorphism classes of planar trees and planar trees with a marked vertex respectively. The middle map forgets the marked point, the source map is given by the number of incoming edges for the marked point and the target map associates to a tree its number of leaves.   			
	The multiplication in this monad is generated by insertion of a tree inside a marked point. 
\end{example}

Let $\mathcal{E}$ be a cocomplete symmetric monoidal category and $T$ be a polynomial functor. One can construct a functor $T^{\mathcal{E}}:\mathcal{E}^J \to \mathcal{E}^I$ given by a formula similar to \ref{PPP}:
\[
T^{\mathcal{E}}(X)_i = \coprod_{b\in t^{-1}(i)} \bigotimes_{e\in p^{-1}(b)} X_{s(e)}.
\]
If $I = J$ and $T$ was given the structure of a polynomial monad then $T^{\mathcal{E}}$ would acquire a structure of monad on $\mathcal{E}^I$.

\begin{definition}
	Let $\mathcal{E}$ be a cocomplete symmetric monoidal category and $T$ a polynomial monad. The category of algebras of $T$ in $\mathcal{E}$ is the category of algebras of the monad $T^{\mathcal{E}}$.
\end{definition}

If $\mathcal{E}$ is the category $\Cat$ of small categories (resp. the category $\SSet$ of simplicial sets), we will call an algebra of $T$ in $\mathcal{E}$ a \emph{categorical} (resp. \emph{simplicial}) $T$-algebra respectively. Explicitly, an algebra $A$ in $\mathcal{E}$ of a polynomial monad $T$ is given by a collection $A_i \in \mathcal{E}$, $i \in I$, equipped with structure maps:
\begin{equation}\label{algebrapolynomialmonad}
m_b: \bigotimes_{e \in p^{-1}(b)} A_{s(e)} \to A_{t(b)}
\end{equation}
for each $b \in B$, which satisfy associativity and unity axioms \cite{bataninberger}.

\subsection{Internal algebras classifiers}\label{sectioninternalalg}

The category of categorical $T$-algebras, which we will write $\Alg_T(\Cat)$, is naturally a $2$-category since it is isomorphic to the category of internal categories in the category of $T$-algebras in $\Set$. The terminal internal category has a unique structure of $T$-algebra in $\Set$ for any polynomial monad $T$; the latter promotes it to a terminal categorical $T$-algebra. Given a morphism of polynomial monads $f: S \to T$ we have a restriction $2$-functor $f^*: \Alg_T(\Cat)\to \Alg_S(\Cat)$. Let us recall the notions of internal algebras \cite{batanin,bataninberger}:
%
%
%
\begin{definition}\label{intalg2}
	Let $f: S \to T$ be a morphism of polynomial monads and $A$ be a categorical $T$-algebra. An \emph{internal $S$-algebra in $A$} is a {lax morphism} of categorical $S$-algebras from the {terminal} categorical $S$-algebra to $f^*(A)$. Internal $S$-algebras in $A$ and $S$-natural transformations form a category $\Int_S(A)$ and this construction extends to a $2$-functor 
	\begin{equation}\label{Int}
	\Int_S:\Alg_T(\Cat)\to\Cat.
	\end{equation}
\end{definition}

For any morphism of polynomial monads $f: S\to T$ one can associate a  categorical $T$-algebra with certain  universal property \cite{batanin,bataninberger,weber}. Namely,
this is the object representing the $2$-functor
\ref{Int}. This categorical $T$-algebra is called \emph{the classifier of internal $S$-algebras inside categorical $T$-algebras} and is denoted $T^S$.
In particular, if $f= Id$ the $T$-algebra $T^T$ is called \emph{absolute classifier of $T$}.

\subsection{Homotopically cofinal maps of polynomial monads}

For a polynomial monad $T$, the category of simplicial $T$-algebras is a simplicial model category, with the model structure transferred from the projective model structure on the category of collections of simplicial sets \cite{bataninberger}. Also note that, if $A$ is a categorical $T$-algebra, then we get a simplicial $T$-algebra, which we will write $N(A)$, by taking the nerve of each category of the underlying collection. We have the following lemma \cite[Corollary 5.2]{batanindeleger}.

\lem\label{lemmaformulaleftderivedfunctor}
Let $f: S \to T$ be a morphism of polynomial monads. Let $f^*$ be the restriction functor from simplicial $T$-algebras to simplicial $S$-algebras and $f_!$ its left adjoint. There is a weak equivalence
\[
\mathbb{L} f_! (1) \sim N \left( T^S \right).
\]
where $\mathbb{L} f_!$ is the left derived functor of $f_!$.
\endlem

This motivates the following definition, which is central in the homotopy theory for polynomial monads \cite{batanindeleger}.

\begin{definition}\label{definitionhomotopicallycofinal}
	A morphism of polynomial monads $f: S\to T$ is called \emph{homotopically cofinal} if $N\left(T^S\right)$ is contractible.
\end{definition}

\subsection{Computation of classifiers}

Let $f: S \to T$ be a morphism of polynomial monad. We have the following commutative square of adjunctions:
\[
\xymatrix{
	\Alg_S \ar@<-.5ex>[rr]_-{f_!} \ar@<-.5ex>[dd]_-{\scriptstyle \mathcal{U}_S} && \Alg_T \ar@<-.5ex>[ll]_-{f^*} \ar@<-.5ex>[dd]_-{\scriptstyle \mathcal{U}_T} \\
	\\
	\Set^J \ar@<-.5ex>[rr]_-{\phi_!} \ar@<-.5ex>[uu]_-{\scriptstyle \mathcal{F}_S} && \Set^I \ar@<-.5ex>[ll]_-{\phi^*} \ar@<-.5ex>[uu]_-{\scriptstyle \mathcal{F}_T}
}
\]
where $\phi^*$ is the  restriction functor induced by $\phi: J \to I$ and $\phi_!: \Set^J\to \Set^I$ is its left adjoint given by coproducts over fibres of $\phi$. It was proved in \cite{batanin,bataninberger} that the classifier induced by $f$ is given by an internal categorical object in the category of $T$-algebras in $\Set$:
\begin{equation}\label{equationclassifier}
\xymatrix{
	\F_T(\phi_!(1)) \ar[r]|-i & \F_T(\phi_!(S 1)) \ar@<-1ex>[l]_-s \ar@<1ex>[l]^-t & \F_T(\phi_!(S^2 1)) \ar[l]|-m \ar@<-1ex>[l]_-{p_2} \ar@<1ex>[l]^-{p_1}
}
\end{equation}
where $1$ is the terminal $J$-collection of sets. The source map $s$ is given using the map $f$ and the multiplication of the monad $T$. The target map $t$ is given by applying $\F_T \cdot \phi_!$ to the unique map. The other structure maps are described in \cite{batanindeleger}.

\begin{remark}\label{remarkcomputationclassifier}
	When the map $\phi$ is the identity, the classifier is particularly easy to compute. As it was explained \cite{bataninberger}, the formula \ref{equationclassifier} tells us that the objects of $(T^S)_i$ are the operations $b$ of $T$ such that $t(b) = i$. The morphisms are the operations $b$ of $T$ such that $t(b) = i$ together with operations $c_1,\ldots,c_k$ of $S$, satisfying $t(c_i) = s(\sigma(i))$ for all $i = 1,\ldots,k$, where $\sigma: \{1,\ldots,k\} \xrightarrow{\sim} p^{-1}(b)$ is an ordering of the fibre. The source of such morphism is obtained by composing $b$ and $\psi(c_1),\ldots,\psi(c_k)$, where $\psi$ is the map as in \ref{morphismpolymon}. The target is $b$.
\end{remark}

\section{Cofinality result}\label{sectioncofinality}

\subsection{Pointed infinitesimal bimodules over operadic bimodules}\label{subsectionpointedinfinitesimalbimodule}

Recall that, if $A$ and $B$ are two non-symmetric operads in a symmetric monoidal category $(\mathcal{E},\otimes,I)$, an $A-B$-bimodule $C$ is given by
\begin{itemize}
	\item a collection $C_n$ of objects  in $\mathcal{E}$ for all integers $n \geq 0$,
	
	\item for $n \geq 0$, $m_1,\ldots,m_n \geq 0$ and $m := m_1 + \ldots m_n$, maps
	\begin{equation}\label{mapbimoduleleft}
	\lambda: A_n \otimes \bigotimes_{j=1}^n C_{m_j} \to C_m
	\end{equation}
	and
	\begin{equation}\label{mapbimoduleright}
	\rho: C_n \otimes \bigotimes_{j=1}^n B_{m_j} \to C_m
	\end{equation}
	such that the usual associativity and unity conditions are satisfied.
\end{itemize}


\begin{definition}\label{definitionweakbimodule}
	Let $A$ and $B$ be two non-symmetric operads in a symmetric monoidal category $(\mathcal{E},\otimes,I)$. Let $C$ and $D$ be two $A-B$-bimodules. An infinitesimal $C-D$-bimodule $E$ is given by
	\begin{itemize}
		\item a collection of objects $E_n$ in $\mathcal{E}$ for $n \geq 0$
		
		\item for any $n \geq 0$, $m_1,\ldots,m_n \geq 0$ with $m := m_1 + \ldots m_n$, and $k \in \{1,\ldots,n\}$, a map
		\begin{equation}\label{mapinfinitesimalleft}
		\phi : A_n \otimes \bigotimes_{j=1}^n (E^{j,k})_{m_j} \to E_m
		\end{equation}
		where we write
		\[
		E^{j,k} :=
		\left\{
		\begin{array}{ll}
		C & \text{if } j < k\\
		E & \text{if } j = k\\
		D & \text{if } j > k
		\end{array}
		\right.
		\]
		
		\item for any $n \geq 0$, and $m_1,\ldots,m_n \geq 0$ with $m := m_1 + \ldots m_n$, a map
		\begin{equation}\label{mapinfinitesimalright}
		\psi : E_n \otimes \bigotimes_{j=1}^n B_{m_j} \to E_m
		\end{equation}
	\end{itemize}
	satisfying the usual unit and associativity axioms, that is the composites
	\[
		E_n \simeq I \otimes E_n \xrightarrow{\epsilon \otimes id} A_1 \otimes E_n \xrightarrow{\phi} E_n
	\]
	and
	\[
		E_n \simeq E_n \otimes \bigotimes_{i=1}^n I \xrightarrow{id \otimes \bigotimes_{i=1}^n \epsilon} E_n \otimes \bigotimes_{i=1}^n B_1 \xrightarrow{\psi} E_n
	\]
	are the identity on $E_n$ and for any $p \geq 0$, $n_1,\ldots,n_p \geq 0$ with $n := n_1+\ldots+n_p$, and $m_1,\ldots,m_n \geq 0$ with $m := m_1+\ldots+m_n$, $m'_i:=m_{n'_{i-1}+1}+\ldots+m_{n'_i}$ and $n'_i:=n_1+\ldots+n_i$ for $i=1,\ldots,p$, $n'_0:=0$,
	\begin{itemize}
		\item if $1 \leq l \leq p$ and $n'_{l-1}+1 \leq k \leq n'_l$, the following diagram commutes
		\[
		\xymatrix{
			A_p \otimes \bigotimes_{i=1}^p A_{n_i} \otimes \bigotimes_{j=1}^n (E^{j,k})_{m_j} \ar[rr]^-\simeq \ar[d]_{\mu \otimes 1} && A_p \otimes \bigotimes_{i=1}^p \left(A_{n_i} \otimes \bigotimes_{j=n'_{i-1}+1}^{n'_i} (E^{j,k})_{m_j} \right) \ar[d]^-{1 \otimes \lambda \otimes \ldots \otimes \phi \otimes \ldots \otimes \lambda} \\
			A_n \otimes \bigotimes_{j=1}^n (E^{j,k})_{m_j} \ar[rd]_{\phi} && A_p \otimes \bigotimes_{i=1}^p (E^{i,l})_{m'_i} \ar[ld]^-{\phi} \\
			& E_m
		}
		\]
		\item if $1 \leq l \leq p$, the following diagram commutes
		\[
		\xymatrix{
			A_p \otimes \bigotimes_{i=1}^p (E^{i,l})_{n_i} \otimes \bigotimes_{j=1}^n B_{m_j} \ar[rr]^-\simeq \ar[d]_{\phi \otimes 1} && A_p \otimes \bigotimes_{i=1}^p \left((E^{i,l})_{n_i} \otimes \bigotimes_{j=n'_{i-1}+1}^{n'_i} B_{m_j} \right) \ar[d]^-{1 \otimes \rho \otimes \ldots \otimes \psi \otimes \ldots \rho} \\
			E_n \otimes \bigotimes_{j=1}^n B_{m_j} \ar[rd]_{\psi} && A_p \otimes \bigotimes_{i=1}^p (E^{i,l})_{m'_i} \ar[ld]^-{\phi} \\
			& E_m
		}
		\]
		\setstretch{1}
		
		\item the following diagram commutes
		\[
		\xymatrix{
			E_p \otimes \bigotimes_{i=1}^p B_{n_i} \otimes \bigotimes_{j=1}^n B_{m_j} \ar[rr]^-\simeq \ar[d]_{\psi \otimes 1} && E_p \otimes \bigotimes_{i=1}^p \left(B_{n_i} \otimes \bigotimes_{j=n'_{i-1}+1}^{n'_i} B_{m_j} \right) \ar[d]^-{1 \otimes \mu \otimes \ldots \otimes \mu} \\
			E_n \otimes \bigotimes_{j=1}^n B_{m_j} \ar[rd]_{\psi} && E_p \otimes \bigotimes_{i=1}^p B_{m'_i} \ar[ld]^-{\psi} \\
			& E_m
		}
		\]
	\end{itemize}
\end{definition}

Note that there is a classical notion of {infinitesimal bimodule} $E$ over a non-symmetric operad $A$ in $\mathcal{E}$. It is given by
\begin{itemize}
	\item a collection of objects $E_n$ in $\mathcal{E}$ for $n \geq 0$
	
	\item for $1 \leq k \leq n$ and $m \geq 0$, maps
	\[
	\bullet_k : A_n \otimes E_m \to E_{m+n-1}
	\]
	
	\item for $1 \leq k \leq n$ and $m \geq 0$, maps
	\[
	\star_k : E_n \otimes A_m \to E_{m+n-1}
	\]
\end{itemize}
satisfying unity, associativity and compatibility axioms \cite{turchin}.

\lem\label{theoremequivalenceweakbimodules}
Let $A$ be a non-symmetric operad. $A$ is a bimodule over itself and the category of infinitesimal $A-A$-bimodules in the sense of Definition \ref{definitionweakbimodule} is equivalent to the category of infinitesimal bimodules over $A$ in the classical sense.
\endlem

\pf
The proof is inspired from the proof of equivalence between May operads and Markl operads \cite{markl2008operads}. First assume that $E$ is an infinitesimal bimodule over $A$ in the classical sense. Let $n \geq 0$, $m_1,\ldots,m_n \geq 0$ and $m = m_1+\ldots+m_n$
%
%
%
We define $\psi$ as the composite
\[
E_n \otimes \bigotimes_{j=1}^n A_{m_j} \xrightarrow{\star_1 \otimes 1} E_{n+m_1-1} \otimes \bigotimes_{j=2}^n A_{m_j} \to \ldots \to E_{m-m_n+1} \otimes A_{m_n} \xrightarrow{\star_{m-m_n+1}} E_m.
\]
and, for $k \in \{1,\ldots,n\}$, $\phi$ as the composite
\[
A_n \otimes \bigotimes_{j=1}^n (E^{j,k})_{m_j} \simeq A_n \otimes \bigotimes_{j \in [n]-k} A_{m_j} \otimes E_{m_k} \xrightarrow{\kappa \otimes 1} A_{m-m_k+1} \otimes E_{m_k} \xrightarrow{\bullet_{m_1+\ldots+m_{k-1}+1}} E_m,
\]
where $[n]-k = \{1,\ldots,n\} \backslash \{k\}$ and $\kappa$ is the composite
\[
A_n \otimes \bigotimes_{j \in [n]-k} A_{m_j} \xrightarrow{1 \otimes \epsilon \otimes 1} A_n \otimes \bigotimes_{j=1}^{k-1} A_{m_j} \otimes A_1 \otimes \bigotimes_{j=k+1}^n A_{m_j} \xrightarrow{\mu} A_{m-m_k+1}.
\]

We get a functor from the category of infinitesimal bimodules over $A$ to the category of infinitesimal $A-A$-bimodules.

In the other direction, assume that $E$ is an infinitesimal $A-A$-bimodule in the sense of Definition \ref{definitionweakbimodule}. Let $1 \leq k \leq n$ and $m \geq 0$. Let $m_k = m$ and $m_j=1$, for $j = 1,\ldots,n$ with $j \neq k$. We define $\bullet_k$ as the composite
\[
A_n \otimes E_m \xrightarrow{1 \otimes \epsilon^{\otimes k-1} \otimes 1 \otimes \epsilon^{\otimes n-k}} A_n \otimes \bigotimes_{j=1}^n (E^{j,k})_{m_j} \xrightarrow{\phi} E_{m+n-1},
\]
and, similarly, $\star_k$ as the composite
\[
E_n \otimes A_m \xrightarrow{1 \otimes \epsilon^{\otimes k-1} \otimes 1 \otimes \epsilon^{\otimes n-k}} E_n \otimes \bigotimes_{j=1}^n A_{m_j} \xrightarrow{\psi} E_{m+n-1}.
\]

We get a functor from the category of infinitesimal $A-A$-bimodules to the category of infinitesimal bimodules over $A$.

It is easy to check that the two functors we have constructed are inverse of each other. This concludes the proof.
\epf

In the rest of this subsection, we will give a notion of pointedness for bimodules and infinitesimal bimodules.

\begin{definition}\label{definitionpointedbimodules}
	Let $A$ and $B$ be two non-symmetric operads in $\mathcal{E}$, and $C$ an $A-B$-bimodule. We say that $C$ is \emph{pointed} if it is equipped with a map $A \to C$ of left $A$-modules and a map $B \to C$ of right $B$-modules such that the following diagram commutes
	\begin{equation}\label{axiompointedbimodule}
	\xymatrix{
		I \ar[r] \ar[d] & A_1 \ar[d] \\
		B_1 \ar[r] & C_1
	}
	\end{equation}
\end{definition}

Note that for $A$ and $B$ two non-symmetric operads in a symmetric monoidal category $\mathcal{E}$ and an $A-B$-bimodule $C$, there is an obvious notion of infinitesimal left $C$-module. It is given by a collection of objects $E_n$ in $\mathcal{E}$ for $n \geq 0$, together with maps as in \ref{mapinfinitesimalleft} with $k=n$ and maps as in \ref{mapinfinitesimalright}, satisfying the same axioms as in Definition \ref{definitionweakbimodule}. Similarly, one has a notion of infinitesimal right $C$-module. Also note that $C$ is an infinitesimal $C-C$-bimodule, and so in particular, an infinitesimal left and right $C$-module.

\begin{definition}\label{definitionpointedweakbimodule}
	Let $A$ and $B$ be two non-symmetric operads and $C$ and $D$ two pointed $A-B$-bimodules. An infinitesimal $C-D$-bimodule $E$ is \emph{pointed} if it is equipped with a map $C \to E$ of infinitesimal left $C$-modules and a map $D \to E$ of infinitesimal right $D$-modules such that, for all $n \geq 0$, $m_1,\ldots,m_n \geq 0$, $m=m_1+\ldots+m_n$ and $k = 1,\ldots,n-1$, the following diagram commutes
	\[
	\xymatrix{
		A_n \otimes \bigotimes_{j=1}^k C_{m_j} \otimes \bigotimes_{j=k+1}^n D_{m_j} \ar[r] \ar[d] & A_n \otimes \bigotimes_{j=1}^n (E^{j,k})_{m_j} \ar[d] \\
		A_n \otimes \bigotimes_{j=1}^n (E^{j,k+1})_{m_j} \ar[r] & E_m
	}
	\]
	and the following diagram commutes in the category of collections
	\begin{equation}\label{mapspointedness}
	\xymatrix{
		& B \ar[rd] \ar[ld] \\
		C \ar[r] & E & D \ar[l] \\
		& A \ar[ru] \ar[lu]
	}
	\end{equation}
	where the maps from $A$ and $B$ are given by Definition \ref{definitionpointedbimodules}.
\end{definition}

\subsection{Description of the morphism of polynomial monads}\label{subsectiondescriptionpolynomialmonads}

Let $\wbimodbfive$ be the polynomial monad for quintuples $(A,B,C,D,E)$, where $A$ and $B$ are non-symmetric operads, $C$ and $D$ are pointed $A-B$-bimodules and $E$ is a pointed infinitesimal $C-D$-bimodule, and $\nopbfive$ be the polynomial monad for commutative diagrams
\begin{equation}\label{diamonddiagram}
\xymatrix{
	& B \ar[rd] \ar[d] \ar[ld] \\
	C \ar[r] & E & D \ar[l] \\
	& A \ar[ru] \ar[u] \ar[lu]
}
\end{equation}
of non-symmetric operads.

The objective of this section is to prove the following theorem.
\[
\xymatrix{
	\Alg_S(\Set) \ar@<-.5ex>[rr]_-{f_*} \ar@<-.5ex>[dd]_-{\scriptstyle \mathcal{U}_S} && \Alg_T(\Set) \ar@<-.5ex>[ll]_-{f^*} \ar@<-.5ex>[dd]_-{\scriptstyle \mathcal{U}_T} \\
	\\
	\Set^J \ar@<-.5ex>[rr]_-{\phi_*} \ar@<-.5ex>[uu]_-{\scriptstyle \mathcal{F}_S} && \Set^I \ar@<-.5ex>[ll]_-{\phi^*} \ar@<-.5ex>[uu]_-{\scriptstyle \mathcal{F}_T}
}
\]

\thm\label{theoremcofinality}
There is a homotopically cofinal morphism of polynomial monads
\begin{equation}\label{homotopycofinalmap}
f: \wbimodbfive \to \nopbfive
\end{equation}
\endthm

First, let us give an explicit description of $\nopbfive$. In order to do this, we define a partial order on the set $\{A,B,C,D,E\}$.

\begin{definition}\label{abcderelation}
	For any $l_1,l_2 \in \{A,B,C,D,E\}$, we write $l_1 \leq l_2$ if there is an arrow from $l_1$ to $l_2$ in the diagram \ref{diamonddiagram}. Explicitly, we have the following relations:
	\[
	\begin{array}{cc}
	A \leq C \leq E & B \leq C \leq E \\
	A \leq D \leq E & B \leq D \leq E
	\end{array}
	\]
\end{definition}

Let $\mathbb{N}$ be the set of non-negative integers. The polynomial monad $\nopbfive$ is given by the polynomial
\[
\xymatrix{
	\{A,B,C,D,E\} \times \mathbb{N} & {PTr}_5^* \ar[r] \ar[l] & PTr_5 \ar[r] & \{A,B,C,D,E\} \times \mathbb{N}
}
\]
where $PTr_5$ is the set of isomorphism classes of planar trees where each vertex $v$ has a label $l_v \in \{A,B,C,D,E\}$. Each tree itself also has a label $l \in \{A,B,C,D,E\}$, called \emph{target label}, satisfying $l_v \leq l$ for each vertex $v$ of the tree. $PTr_5^*$ is the set of elements of $PTr_5$ with one vertex marked. The source map produces a label corresponding to the label of the marked vertex, and a number equal to the number of incoming edges of this vertex. The target map returns the target label and the number of leaves. The multiplication is by insertion of trees inside the vertices of another tree. The condition of insertion is that the number of leaves and the target label of each inserted tree must correspond to the number of incoming edges and the label of the vertex where it is inserted. We will not prove that the polynomial monad just described gives the correct category of algebras here, because the proof would be very similar to the proof of Proposition \ref{propositionpolymongivescorrectalgebras}.

%

Now we will describe in more detail the polynomial monad $\wbimodbfive$.

\begin{definition}\label{conditionstar}
	We say that the labels $l_v \in \{A,B,C,D,E\}$ of the vertices $v$ of an isomorphism class of planar trees \emph{lie on a line and a point} if, for a representant of the isomorphism class, it is possible to draw a line and a point on this line such that
	\begin{itemize}
		\item all the vertices below the line have label $A$
		\item all the vertices above the line have label $B$
		\item all the vertices on the line and to the left of the point have label $C$
		\item all the vertices on the line and to the right of the point have label $D$
		\item if a vertex lies on the point, it has label $E$
	\end{itemize}
	
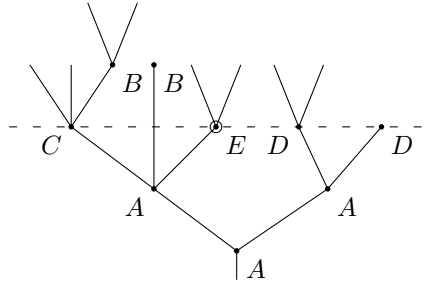
	\[
	\begin{tikzpicture}[scale=.55]
	\draw[loosely dashed] (-5,0) -- (5,0);
	\draw[fill=white] (0,0) circle (4pt);
	
	\draw (-3.5,0) -- (-1.5,-1.5);
	\draw (0,0) -- (-1.5,-1.5);
	\draw (-1.5,-1.5) -- (.5,-3);
	\draw (2,0) -- (2.7,-1.5);
	\draw (4,0) -- (2.7,-1.5);
	\draw (.5,-3) -- (.5,-3.7);
	\draw (-1.5,-1.5) -- (-1.5,1.5);
	\draw (2.7,-1.5) -- (.5,-3);
	
	\draw (-3.5,0) -- (-4.5,1.5);
	\draw (-3.5,0) -- (-3.5,1.5);
	\draw (-3.5,0) -- (-2.5,1.5);
	\draw (-2.5,1.5) -- (-3.1,3);
	\draw (-2.5,1.5) -- (-1.9,3);
	
	\draw (0,0) -- (-.6,1.5);
	\draw (0,0) -- (.6,1.5);
	
	\draw (2,0) -- (1.4,1.5);
	\draw (2,0) -- (2.6,1.5);
	
	\draw[fill] (.5,-3) circle (1.5pt) node[below right]{$A$};
	\draw[fill] (-1.5,-1.5) circle (1.5pt) node[below left]{$A$};
	\draw[fill] (2.7,-1.5) circle (1.5pt) node[below right]{$A$};
	\draw[fill] (-3.5,0) circle (1.5pt) node[below left]{$C$};
	\draw[fill] (0,0) circle (1.5pt) node[below right]{$E$};
	\draw[fill] (2,0) circle (1.5pt) node[below left]{$D$};
	\draw[fill] (4,0) circle (1.5pt) node[below right]{$D$};
	\draw[fill] (-2.5,1.5) circle (1.5pt) node[below right]{$B$};
	\draw[fill] (-1.5,1.5) circle (1.5pt) node[below right]{$B$};
	
	\end{tikzpicture}
	\]
\end{definition}

The polynomial monad $\wbimodbfive$ is given by the polynomial
\[
\xymatrix{
	\{A,B,C,D,E\} \times \mathbb{N} & IPTr_5^* \ar[r] \ar[l] & IPTr_5 \ar[r] & \{A,B,C,D,E\} \times \mathbb{N}
}
\]
where $IPTr_5$ is the subset of $PTr_5$ of isomorphisms classes of planar trees whose labels lie on a line and a point. The rest of the description of the polynomial monad $\wbimodbfive$ is straightforward, as it is completely similar to the description of the polynomial monad $\nopbfive$.

\begin{proposition}\label{propositionpolymongivescorrectalgebras}
	The polynomial monad $\wbimodbfive$ described above gives the desired category of algebras.
\end{proposition}

\pf
	An algebra over $\wbimodbfive$ is given by five collections over $\mathbb{N}$, which we will denote by $A$, $B$, $C$, $D$ and $E$, equipped with structure maps as in \ref{algebrapolynomialmonad}. The maps \ref{mapnonsymmetricoperad}, \ref{mapbimoduleleft}, \ref{mapbimoduleright}, \ref{mapinfinitesimalleft} and \ref{mapinfinitesimalright} correspond to structure maps induced by some of the two levelled trees in $IPTr_5$. For example, the maps \ref{mapinfinitesimalleft} correspond to the structure maps induced by two levelled trees whose vertex at the bottom has label $A$, a unique vertex at the top has label $E$, all vertices to the left of it have label $C$, all vertices to the right of it have label $D$, as in the following picture:
	\[
	\begin{tikzpicture}[scale=.8]
	\draw (0,-.5) -- (0,0);
	\draw (-1.8,1) -- (0,0) -- (1.8,1);
	\draw (-.6,1) -- (0,0) -- (.6,1) -- (.6,1.7);
	
	\draw (-2.1,1.7) -- (-1.8,1) -- (-1.5,1.7);
	\draw (1.1,1.7) -- (.6,1) -- (.1,1.7);
	\draw (1.8,1) -- (1.8,1.7);
	
	\draw[fill] (0,0) circle (1pt) node[below right]{$A$};
	\draw[fill] (-1.8,1) circle (1pt) node[left]{$C$};
	\draw[fill] (-.6,1) circle (1pt) node[left]{$C$};
	\draw[fill] (.6,1) circle (1pt) node[right]{$E$};
	\draw[fill] (1.8,1) circle (.8pt) node[right]{$D$};
	\end{tikzpicture}
	\]
	The maps coming from pointedness, as in \ref{mapspointedness}, correspond to structure maps induced by the corollas in $IPTr_5$. Note that for a pointed $A-B$-bimodule $C$, the following square commutes:
	\[
	\xymatrix{
		A_n \otimes \bigotimes_{j=1}^n B_{m_j} \ar[r] \ar[d] & A_n \otimes \bigotimes_{j=1}^n C_{m_j} \ar[d] \\
		C_n \otimes \bigotimes_{j=1}^n B_{m_j} \ar[r] & C_m
	}
	\]
	This, together with the axioms of pointed infinitesimal bimodules of Definition \ref{definitionpointedbimodules}, ensures us that the structure maps induced by two levelled trees in $IPTr_5$ correspond without ambiguity to a map obtained by combining maps coming from the pointed infinitesimal bimodule structure. The associativity axioms ensure us that the structure maps induced by any tree in $IPTr_5$ correspond without ambiguity to a map obtained by combining maps coming from the pointed infinitesimal bimodule structure.
\epf

The map $f: \wbimodbfive \to \nopbfive$ is given by a diagram as \ref{morphismpolymon}, where the vertical maps are inclusions of sets.

\subsection{Description of the induced classifier}\label{subsectionclassifier}

Let us describe the classifier associated to the map \ref{homotopycofinalmap}. It is given by a collection of categories indexed by the set of colours of $\nopbfive$, that is $\{A,B,C,D,E\} \times \mathbb{N}$. Let $(l,n)$ be an element of this set. The category indexed by $(l,n)$ has the following description.

According to Remark \ref{remarkcomputationclassifier}, the objects are elements of $PTr_5$ whose target is $(l,n)$, that is isomorphism classes of planar tree with $n$ leaves and target label $l$, where each vertex $v$ has a label $l_v \in \{A,B,C,D,E\}$ satisfying $l_v \leq l$. We will describe morphisms as nested trees. It is an isomorphism class of planar tree $T$ where each vertex itself contains an isomorphism class of planar tree. The number of incoming edges of each vertex must be equal to the number of leaves of the tree inside it. The tree $T$ has $n$ leaves and each vertex $v$ has a label $l_v \in \{A,B,C,D,E\}$ satisfying $l_v \leq l$. Each vertex $w$ of each tree inside a vertex $v$ also has a label $l_w$ such that $l_w \leq l_v$. Finally, the vertices of a tree inside a vertex must satisfy the condition of Definition \ref{conditionstar}.
Here is an example of nested tree:
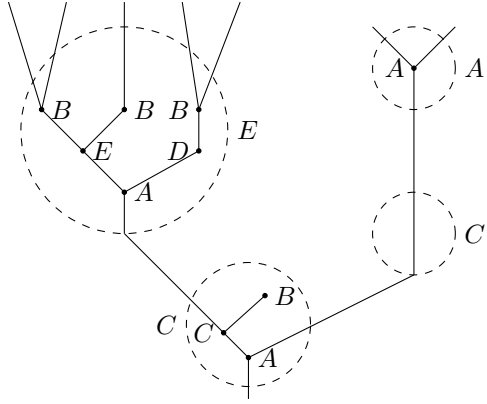
\[
\begin{tikzpicture}[scale = .55]
\draw (0,-1) -- (0,0);

\draw[fill] (0,0) circle (1.5pt) node[right]{$A$};

\draw (0,0) -- (-3,3);

\draw[dashed] (-3,5.5) circle (2.5);
\draw (-.5,5.5) node[right]{$E$};

\draw (-3,3) -- (-3,4);
\begin{scope}[shift={(0,-.8)}]
\draw (-3,4.8) -- (-4,5.8);
\draw (-3,4.8) -- (-1.2,5.8);
\begin{scope}[shift={(.4,0)}]
\draw (-4.4,5.8) -- (-5.4,6.8);
\draw (-4.4,5.8) -- (-3.4,6.8);
\draw (-1.6,5.8) -- (-1.6,6.8);
\draw (-5.4,6.8) -- (-6.2,9.4);
\draw (-5.4,6.8) -- (-4.8,9.4);
\draw (-3.4,6.8) -- (-3.4,9.4);
\draw (-1.6,6.8) -- (-2,9.4);
\draw (-1.6,6.8) -- (-.6,9.4);

\draw[fill] (-3.4,4.8) circle (1.5pt) node[right]{$A$};
\draw[fill] (-4.4,5.8) circle (1.5pt) node[right]{$E$};
\draw[fill] (-1.6,5.8) circle (1.5pt) node[left]{$D$};
\draw[fill] (-5.4,6.8) circle (1.5pt) node[right]{$B$};
\draw[fill] (-3.4,6.8) circle (1.5pt) node[right]{$B$};
\draw[fill] (-1.6,6.8) circle (1.5pt) node[left]{$B$};
\end{scope}
\end{scope}

\draw (0,0) -- (4,2);

\draw (4,2) -- (4,7);
\draw (4,7) -- (3,8);
\draw (4,7) -- (5,8);

\draw[fill] (4,7) circle (1.5pt) node[left]{$A$};
\draw[dashed] (4,7) circle(1);
\draw[dashed] (4,3) circle(1);
\draw (5,3) node[right]{$C$};
\draw (5,7) node[right]{$A$};

\draw (-.6,.6) -- (.4,1.5);
\draw[fill] (.4,1.5) circle (1.5pt) node[right]{$B$};
\draw[fill] (-.6,.6) circle (1.5pt) node[left]{$C$};
\draw[dashed] (0,.8) circle (1.5);
\draw (-1.5,.8) node[left]{$C$};

\end{tikzpicture}
\]
The tree $T$ corresponds to $b$ in Remark \ref{remarkcomputationclassifier}, while the trees inside each vertex of $T$ corresponds to $c_1,\ldots,c_k$. The source of a nested tree is obtained by inserting each tree into the vertex it decorates and the target is obtained by forgetting the trees inside each vertex. For example, the nested tree of the previous picture represents the following morphism:

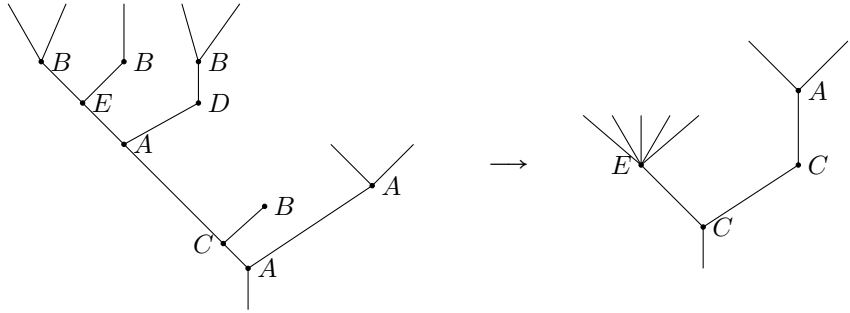
\[
\begin{tikzpicture}[scale = .55]
\draw (0,-1) -- (0,0);

\draw[fill] (0,0) circle (1.5pt) node[right]{$A$};

\draw (0,0) -- (-3,3);

\begin{scope}[shift={(0,-1.8)}]
\draw (-3,4.8) -- (-4,5.8);
\draw (-3,4.8) -- (-1.2,5.8);
\begin{scope}[shift={(.4,0)}]
\draw (-4.4,5.8) -- (-5.4,6.8);
\draw (-4.4,5.8) -- (-3.4,6.8);
\draw (-1.6,5.8) -- (-1.6,6.8);
\draw (-5.4,6.8) -- (-6.2,8.2);
\draw (-5.4,6.8) -- (-4.8,8.2);
\draw (-3.4,6.8) -- (-3.4,8.2);
\draw (-1.6,6.8) -- (-2,8.2);
\draw (-1.6,6.8) -- (-.6,8.2);

\draw[fill] (-3.4,4.8) circle (1.5pt) node[right]{$A$};
\draw[fill] (-4.4,5.8) circle (1.5pt) node[right]{$E$};
\draw[fill] (-1.6,5.8) circle (1.5pt) node[right]{$D$};
\draw[fill] (-5.4,6.8) circle (1.5pt) node[right]{$B$};
\draw[fill] (-3.4,6.8) circle (1.5pt) node[right]{$B$};
\draw[fill] (-1.6,6.8) circle (1.5pt) node[right]{$B$};
\end{scope}
\end{scope}

\draw (0,0) -- (3,2);

\draw (3,2) -- (2,3);
\draw (3,2) -- (4,3);

\draw[fill] (3,2) circle (1.5pt) node[right]{$A$};

\draw (-.6,.6) -- (.4,1.5);
\draw[fill] (.4,1.5) circle (1.5pt) node[right]{$B$};
\draw[fill] (-.6,.6) circle (1.5pt) node[left]{$C$};

\draw (6.3,2.5) node{$\longrightarrow$};

\begin{scope}[shift={(11,1)}]
\draw (0,-1) -- (0,0);

\draw[fill] (0,0) circle (1.5pt) node[right]{$C$};

\draw (0,0) -- (-1.5,1.5);

\draw[fill] (-1.5,1.5) circle (1.5pt) node[left]{$E$};

\draw (-1.5,1.5) -- (-2.9,2.7);
\draw (-1.5,1.5) -- (-2.2,2.7);
\draw (-1.5,1.5) -- (-1.5,2.7);
\draw (-1.5,1.5) -- (-.8,2.7);
\draw (-1.5,1.5) -- (-.1,2.7);

\draw (0,0) -- (2.3,1.5);

\draw (2.3,1.5) -- (2.3,3.3);
\draw (2.3,3.3) -- (1.1,4.5);
\draw (2.3,3.3) -- (3.5,4.5);

\draw[fill] (2.3,3.3) circle (1.5pt) node[right]{$A$};
\draw[fill] (2.3,1.5) circle (1.5pt) node[right]{$C$};

\end{scope}

\end{tikzpicture}
\]
In other words, morphisms are given by contractions of some subtrees to corollas.

\subsection{Proof of the cofinality result}

In order to prove Theorem \ref{theoremcofinality}, we need the Cisinski lemma about smooth functors. For a functor $F: \mathcal{X} \to \mathcal{Y}$ between categories and $y \in \mathcal{Y}$, we will write $F_y$ for the fibre of $F$ over $y$, that is the full subcategory of objects $x \in \mathcal{X}$ such that $F(x) = y$.

\begin{definition}
	A functor $F : \mathcal{X} \to \mathcal{Y}$ is \emph{smooth} if, for all $y \in \mathcal{Y}$, the canonical functor
	\[
	F_y \to y/F
	\]
	induces a weak equivalence between nerves.
	
	Dually, a functor $F : \mathcal{X} \to \mathcal{Y}$ is \emph{proper} if, for all $y \in \mathcal{Y}$, the canonical functor
	\[
	F_y \to F/y
	\]
	induces a weak equivalence between nerves.
\end{definition}

Let us state the Cisinski lemma \cite[Proposition 5.3.4]{cisinski}.

\lem\label{cisinskilemma}
A functor $F : \mathcal{X} \to \mathcal{Y}$ is smooth if and only if for all maps $f_1 : y_0 \to y_1$ in $\mathcal{Y}$ and objects $x_1$ in $\mathcal{X}$ such that $F(x_1) = y_1$, the nerve of the \emph{lifting category} of $f_1$ over $x_1$, whose objects are arrows $f : x \to x_1$ such that $F(f) = f_1$ and morphisms are commutative triangles
\[
\xymatrix{
	x \ar[rr]^g \ar[rd]_f && x' \ar[ld]^{f'} \\
	& x_1
}
\]
with $g$ a morphism in $F_{y_0}$, is contractible.

There is a dual characterisation for proper functors.
\endlem

Now we will construct a functor $F: \mathcal{X} \to \mathcal{Y}$ which we will prove is smooth. We have the following commutative square of polynomial monads
\[
\xymatrix{
	\wbimodbfive \ar[r]^-{u f} \ar[d]_f & \mathrm{NOp} \ar@{=}[d] \\
	\nopbfive \ar[r]_-u & \mathrm{NOp}
}
\]
where $u$ is the morphism of polynomial monads given by the projections. As it was proven in \cite[Proposition 4.7]{batanindeleger}, this commutative square induces a strict map of algebras
\begin{equation}\label{mapbetweenclassifiers}
\nopbfive^{\wbimodbfive} \to u^* \left( \mathrm{NOp}^{\mathrm{NOp}} \right)
\end{equation}
To simplify the notations, let us assume that a colour of $\nopbfive$ is fixed and let $F: \mathcal{X} \to \mathcal{Y}$ be the underlying functor of the map \ref{mapbetweenclassifiers} between the categories indexed by this colour. This functor sends a labelled tree to the same tree, but without labels.

\lem\label{lemmasmoothness}
The functor $F: \mathcal{X} \to \mathcal{Y}$ is smooth.
\endlem

First, let us prove the following technical lemma.

\lem\label{lemmacontractibleliftingcategoryinparticularcase}
Let $f_1: y_0 \to y_1$ be a map in $\mathcal{Y}$ and $x_1$ be an object in $\mathcal{X}$ such that $F(x_1) = y_1$. If $x_1$ is a corolla, that is a tree with exactly one vertex, then the nerve of the lifting category $\mathcal{X}(x_1,f_1)$ of $f_1$ over $x_1$ is contractible.
\endlem

\pf
Let $f_1: y_0 \to y_1$ be a map in $\mathcal{Y}$ and $x_1$ be an object in $\mathcal{X}$ such that $F(x_1) = y_1$. If $x_1$ is a corolla, then $y_1$ is also a corolla and the map $f_1$ is the unique map. Therefore, the pair $(x_1,f_1)$ is uniquely determined by a pair $(l,T)$, where $l \in \{A,B,C,D,E\}$ is the label of the unique vertex of $x_1$ and $T$ is the tree given by the object $y_0 \in \mathcal{Y}$. Let us write $\chi(l,T) = \mathcal{X}(x_1,f_1)$. The objects of $\chi(l,T)$ are the objects $x \in F_{y_0}$ such that there is an arrow $f: x \to x_1$ which satisfies $F(f) = f_1$. Indeed, since the functor $F: \mathcal{X} \to \mathcal{Y}$ is faithful, the map $f$ is automatically given. This means that the objects are labellings of the vertices of $T$ such that the labels lie on a line and a point. The morphisms change the labels as in the diagram \ref{diamonddiagram}. 
For example, if the label $l$ is $E$ and $T$ is the following tree
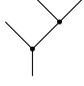
\[
\begin{tikzpicture}
\begin{scope}[scale=.6]
\draw (-.3,-.9) -- (-.3,-.3);
\draw (-.3,-.3) -- (-.9,.3);
\draw (-.3,-.3) -- (.3,.3);
\draw (.3,.3) -- (-.3,.9);
\draw (.3,.3) -- (.9,.9);

\draw[fill] (-.3,-.3) circle (1.2pt);
\draw[fill] (.3,.3) circle (1.2pt);

\end{scope}
\end{tikzpicture}
\]
then $\chi(l,T)$ is represented in the following picture

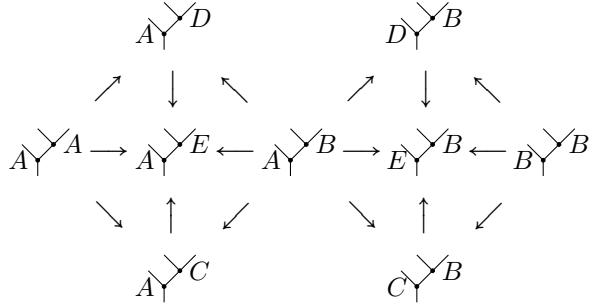
\def\tree#1#2#3{
	\begin{scope}[shift={#1},scale=.5]
		\draw (-.3,-.9) -- (-.3,-.3);
		\draw (-.3,-.3) -- (-.9,.3);
		\draw (-.3,-.3) -- (.3,.3);
		\draw (.3,.3) -- (-.3,.9);
		\draw (.3,.3) -- (.9,.9);

		\ifthenelse{#2=1}
		{
			\draw[fill] (-.3,-.3) node[left]{$\tiny{A}$};
		}{\ifthenelse{#2=2}
			{
				\draw[fill] (-.3,-.3) node[left]{$\tiny{B}$};
				
			}{\ifthenelse{#2=3}
				{
					\draw[fill] (-.3,-.3) node[left]{$\tiny{C}$};
					
				}{\ifthenelse{#2=4}
					{
						\draw[fill] (-.3,-.3) node[left]{$\tiny{D}$};
						
					}{\ifthenelse{#2=5}
						{
							\draw[fill] (-.3,-.3) node[left]{$\tiny{E}$};	
						}{}}}}}

		\ifthenelse{#3=1}
		{
			\draw[fill] (.3,.3) node[right]{$\tiny{A}$};
		}{\ifthenelse{#3=2}
			{
				\draw[fill] (.3,.3) node[right]{$\tiny{B}$};
			}{\ifthenelse{#3=3}
				{
					\draw[fill] (.3,.3) node[right]{$\tiny{C}$};
				}{\ifthenelse{#3=4}
					{
						\draw[fill] (.3,.3) node[right]{$\tiny{D}$};
					}{\ifthenelse{#3=5}
						{
							\draw[fill] (.3,.3) node[right]{$\tiny{E}$};	
						}{}}}}}
		
		\draw[fill] (-.3,-.3) circle (2pt);
		\draw[fill] (.3,.3) circle (2pt);
	\end{scope}
}

\[
\begin{tikzpicture}[scale=.7]

\tree{(-4.8,0)}{1}{1};
\tree{(-2.4,-2.4)}{1}{3};
\tree{(-2.4,0)}{1}{5};
\tree{(-2.4,2.4)}{1}{4};
\tree{(0,0)}{1}{2};
\tree{(2.4,-2.4)}{3}{2};
\tree{(2.4,0)}{5}{2};
\tree{(2.4,2.4)}{4}{2};
\tree{(4.8,0)}{2}{2};

\draw (-3.6,-1.2) node[rotate=-45]{$\longrightarrow$};
\draw (-3.6,0) node{$\longrightarrow$};
\draw (-3.6,1.2) node[rotate=45]{$\longrightarrow$};
\draw (-2.4,-1.2) node[rotate=90]{$\longrightarrow$};
\draw (-2.4,1.2) node[rotate=-90]{$\longrightarrow$};
\draw (-1.2,-1.2) node[rotate=-135]{$\longrightarrow$};
\draw (-1.2,0) node{$\longleftarrow$};
\draw (-1.2,1.2) node[rotate=135]{$\longrightarrow$};
\draw (1.2,-1.2) node[rotate=-45]{$\longrightarrow$};
\draw (1.2,0) node{$\longrightarrow$};
\draw (1.2,1.2) node[rotate=45]{$\longrightarrow$};
\draw (2.4,-1.2) node[rotate=90]{$\longrightarrow$};
\draw (2.4,1.2) node[rotate=-90]{$\longrightarrow$};
\draw (3.6,-1.2) node[rotate=-135]{$\longrightarrow$};
\draw (3.6,0) node{$\longleftarrow$};
\draw (3.6,1.2) node[rotate=135]{$\longrightarrow$};

\end{tikzpicture}
\]

We will prove that the nerve of $\chi(l,T)$ is contractible by induction on the number of vertices of $T$. If $T$ has no vertices, that is the free living edge, then $\chi(l,T)$ is the terminal category and its nerve is contractible. Now assume that $T$ has at least one vertex. Let $v$ be the root vertex and $k$ be the number of incoming edges of $v$. Let $T_1,\ldots,T_k$ be the subtrees above the edges of $v$:

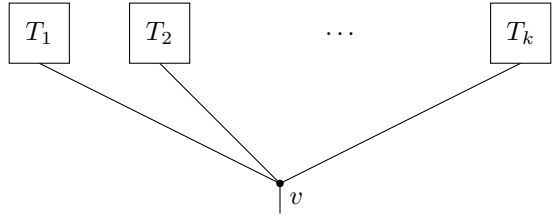
\[
\begin{tikzpicture}[scale=.8]
\draw (0,-.5) -- (0,0);
\draw (0,0) -- (-4,2);
\draw (0,0) -- (-2,2);
\draw (0,0) -- (4,2);
\draw (-4.5,2) rectangle (-3.5,3);
\draw (-2.5,2) rectangle (-1.5,3);
\draw (3.5,2) rectangle (4.5,3);

\draw (-4,2.5) node{$T_1$};
\draw (-2,2.5) node{$T_2$};
\draw (1,2.5) node{$\ldots$};
\draw (4,2.5) node{$T_k$};

\draw[fill] (0,0) circle (1.5pt) node[below right]{$v$};
\end{tikzpicture}
\]

If the label $l$ is $A$ or $B$ then the result is trivial. Otherwise, let $\chi_A(l,T)$ be the full subcategory of $\chi(l,T)$ of trees where the root vertex is labelled with $A$. Let $\chi_B(l,T)$ be the full subcategory of $\chi(l,T)$ of trees where all the vertices other than the root vertex are labelled with $B$. Any simplex of $\chi(l,T)$ is completely contained within $\chi_A(l,T)$ or $\chi_B(l,T)$. Indeed, any object of $\chi(l,T)$ is an object in either $\chi_A(l,T)$ or $\chi_B(l,T)$ and if the last vertex of a simplex is an object of either one of these subcategories, then so must be all the vertices of the simplex. If we can prove that the nerves of $\chi_A(l,T)$, $\chi_B(l,T)$ and the intersection of these two subcategories are contractible, then we can conclude that the nerve of their union is also contracible, which is the desired result. The intersection of $\chi_A(l,T)$ and $\chi_B(l,T)$ is the terminal category. The subcategory $\chi_B(l,T)$ is equivalent to the category $\chi(l,v)$, where $v$ is the corolla with one vertex and $k$ leaves, which has obviously a contractible nerve.

Now, let us consider $\chi_A(l,T)$. If the label $l$ is $C$ or $D$, then this category is equivalent to the product of categories $\prod_{i=1}^k \chi(l,T_i)$, which has contractible nerve by induction. If the label $l$ is $E$, we need to do an induction on $k$.
Let $\chi_C(E,T)$ be the full subcategory of $\chi_A(E,T)$ of trees where the vertices of $T_1$ have label $A$, $B$ or $C$. Let $\chi_D(E,T)$ be the full subcategory of $\chi_A(E,T)$ of trees where the vertices of $T_2,\ldots,T_k$ have label $A$, $B$ or $D$.

Finally, let us write $T'$ for the tree obtained by removing the first edge above the root vertex $v$ and the tree $T_1$. Any simplex of $\chi_A(E,T)$ is completely contained within $\chi_C(E,T)$ or $\chi_D(E,T)$. As before, we want to prove that the nerves of $\chi_C(E,T)$, $\chi_D(E,T)$ and the intersection of these two subcategories are contractible. The subcategory $\chi_C(E,T)$ is equivalent to $\chi(C,T_1) \times \chi_A(E,T')$, the subcategory $\chi_D(E,T)$ is equivalent to $\chi(E,T_1) \times \chi_A(D,T')$ and their intersection is equivalent to $\chi(C,T_1) \times \chi_A(D,T')$. By induction, all three have contractible nerve.
\epf

\pfof[Lemma \ref{lemmasmoothness}]
Let $f_1: y_0 \to y_1$ be a map in $\mathcal{Y}$ and $x_1$ be an object in $\mathcal{X}$ such that $F(x_1) = y_1$. According to the Cisinski lemma \ref{cisinskilemma}, we get the desired result if we can prove that the nerve of the lifting category $\mathcal{X}(x_1,f_1)$ of $f_1$ over $x_1$ is contractible. Remark that $f_1$ can be described by maps $f_1^{(v)} : y_0^{(v)} \to y_1^{(v)}$ for $v \in V$, where $V$ is the set of vertices of $y_1$ and $y_1^{(v)}$ is always the corolla. For example, the morphism
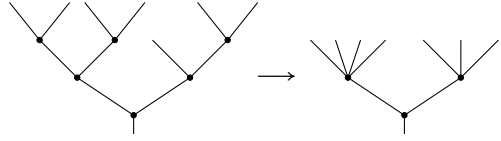
\[
\begin{tikzpicture}[scale=.5]
\draw (0,-.5) -- (0,0);
\draw (0,0) -- (-1.5,1);
\draw (0,0) -- (1.5,1);
\draw (-1.5,1) -- (-2.5,2);
\draw (-1.5,1) -- (-.5,2);
\draw (1.5,1) -- (.5,2);
\draw (1.5,1) -- (2.5,2);
\draw (-2.5,2) -- (-3.3,3);
\draw (-2.5,2) -- (-1.7,3);
\draw (-.5,2) -- (-1.3,3);
\draw (-.5,2) -- (.3,3);
\draw (2.5,2) -- (1.7,3);
\draw (2.5,2) -- (3.3,3);

\draw[fill] (0,0) circle (2pt);
\draw[fill] (-1.5,1) circle (2pt);
\draw[fill] (1.5,1) circle (2pt);
\draw[fill] (-2.5,2) circle (2pt);
\draw[fill] (-.5,2) circle (2pt);
\draw[fill] (2.5,2) circle (2pt);

\draw (3.8,1) node{$\longrightarrow$};

\begin{scope}[shift={(7.2,0)}]
\draw (0,-.5) -- (0,0);
\draw (0,0) -- (-1.5,1);
\draw (0,0) -- (1.5,1);
\draw (-1.5,1) -- (-2.5,2);
\draw (-1.5,1) -- (-1.833,2);
\draw (-1.5,1) -- (-1.166,2);
\draw (-1.5,1) -- (-.5,2);
\draw (1.5,1) -- (.5,2);
\draw (1.5,1) -- (1.5,2);
\draw (1.5,1) -- (2.5,2);

\draw[fill] (0,0) circle (2pt);
\draw[fill] (-1.5,1) circle (2pt);
\draw[fill] (1.5,1) circle (2pt);
\end{scope}
\end{tikzpicture}
\]
can be described by the three morphisms
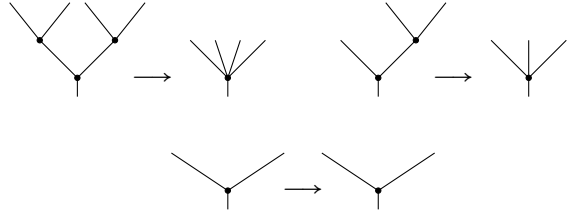
\[
\begin{tikzpicture}[scale=.5]
\draw (0,-.5) -- (0,0);
\draw (0,0) -- (-1.5,1);
\draw (0,0) -- (1.5,1);
\draw[fill] (0,0) circle (2pt);

\draw (2,0) node{$\longrightarrow$};

\begin{scope}[shift={(4,0)}]
\draw (0,-.5) -- (0,0);
\draw (0,0) -- (-1.5,1);
\draw (0,0) -- (1.5,1);
\draw[fill] (0,0) circle (2pt);
\end{scope}

\begin{scope}[shift={(-4,3)}]
\draw (0,-.5) -- (0,0);
\draw (0,0) -- (-1,1);
\draw (0,0) -- (1,1);
\draw (-1,1) -- (-1.8,2);
\draw (-1,1) -- (-.2,2);
\draw (1,1) -- (.2,2);
\draw (1,1) -- (1.8,2);
\draw[fill] (0,0) circle (2pt);
\draw[fill] (-1,1) circle (2pt);
\draw[fill] (1,1) circle (2pt);

\draw (2,0) node{$\longrightarrow$};

\begin{scope}[shift={(4,0)}]
\draw (0,-.5) -- (0,0);
\draw (0,0) -- (-1,1);
\draw (0,0) -- (-.333,1);
\draw (0,0) -- (.333,1);
\draw (0,0) -- (1,1);
\draw[fill] (0,0) circle (2pt);
\end{scope}
\end{scope}

\begin{scope}[shift={(4,3)}]
\draw (0,-.5) -- (0,0);
\draw (0,0) -- (-1,1);
\draw (0,0) -- (1,1);
\draw (1,1) -- (.2,2);
\draw (1,1) -- (1.8,2);
\draw[fill] (0,0) circle (2pt);
\draw[fill] (1,1) circle (2pt);

\draw (2,0) node{$\longrightarrow$};

\begin{scope}[shift={(4,0)}]
\draw (0,-.5) -- (0,0);
\draw (0,0) -- (-1,1);
\draw (0,0) -- (0,1);
\draw (0,0) -- (1,1);
\draw[fill] (0,0) circle (2pt);
\end{scope}
\end{scope}

\end{tikzpicture}
\]
We also have objects $x_1^{(v)}$ for each $v \in V$ which decompose $x_1$, just like $y_1^{(v)}$ decompose $y_1$. The category $\mathcal{X}(x_1,f_1)$ is equivalent to the product of the categories $\mathcal{X}\left(x_1^{(v)},f_1^{(v)}\right)$ over $v \in V$, which has contractible nerve thanks to Lemma \ref{lemmacontractibleliftingcategoryinparticularcase}.
\epf

\pfof[Theorem \ref{theoremcofinality}]
The functor $F: \mathcal{X} \to \mathcal{Y}$ is smooth and its fibres have contractible nerve, since they have a terminal object. Using Quillen Theorem A, we deduce that $F$ induces a weak equivalence between nerve. Again, the nerve of $\mathcal{Y}$ is contractible since this category has a terminal object. So the nerve of $\mathcal{X}$ is contractible, which concludes the proof.
\epf

\section{Turchin-Dwyer-Hess delooping theorem}\label{chapterapplication}

\subsection{Delooping theorem for a left proper model category}\label{subsectionspheresanddisksingeneralmodelcategory}

Let us define, by induction on $n$, spheres $S^n$ and disks $D^n$ in a general model category.

\begin{definition}\label{definitionspheresanddisks}
	Let $M$ be a model category. We will write $S^{-1}$ (resp. $D^{-1}$) for the initial (resp. terminal) object in $M$. Note that there in a unique map $S^{-1} \to D^{-1}$.
	
	Assume we have defined $S^{n-1}$ and $D^{n-1}$ for $n \geq 0$, together with a map $S^{n-1} \to D^{n-1}$. We can factorise this map as follows
	\begin{equation}\label{trianglespheresdisks}
	\xymatrix{
		S^{n-1} \ar[rr] \ar@{ >->}[rd] && D^{n-1} \\
		& D^n \ar[ru]_\sim
	}
	\end{equation}
	where the first map is a cofibration and the second map is a weak equivalence. We define $S^n$ as the pushout
	\begin{equation}\label{definitionsphere}
	\xymatrix{
		S^{n-1} \ar@{ >->}[r] \ar@{ >->}[d] & D^n \ar[d] \\
		D^n \ar[r] & S^n
	}
	\end{equation}
	This comes with a map $S^n \to D^n$ induced by the identity on $D^n$.
\end{definition}

For $n \geq 0$, note that the triangle \ref{trianglespheresdisks} allows us to construct a composite map
\begin{equation}\label{mapgn}
S^{n-1} \rightarrowtail D^n \xrightarrow{\sim} D^{n-1} \xrightarrow{\sim} \ldots \xrightarrow{\sim} D^0
\end{equation}
which we will denote by $g_n$. This induces a Quillen adjunction between comma categories
\[
\xymatrix@C=1pc{
	S^{n-1}/M \ar@/^.7pc/[rr]^{(g_n)_!} & \bot & D^0/M \ar@/^.7pc/[ll]^{g_n^*}
}
\]
where $g_n^*$ is given by precomposition and $(g_n)_!$ by pushout respectively.

We also have a composite
\begin{equation}\label{mapsouthpole}
D^0 \to S^0 \to S^1 \to \ldots \to S^n
\end{equation}
where the first map is one of the two maps to the coproduct and the other maps are given by the commutative square \ref{definitionsphere}. Finally, recall that a model category is \emph{left proper} if weak equivalences are preserved by pushouts along cofibrations.

\lem\label{lemmaweakequivalencespheres}
Let $n \geq 0$ and $M$ be a left proper model category. Then $\mathbb{L}(g_n)_!(1)$, that is the value of the left derived functor of $(g_n)_!$ on the terminal object, is weakly equivalent to the object in $D^0 / M$ given by the map \ref{mapsouthpole}.
\endlem

\pf
We look at the following diagram
\[
\xymatrix{
	S^{n-1} \ar@/^1pc/[rr]^{g_n} \ar@{>->}[r] \ar@{ >->}[d] & D^n \ar[r] \ar@{->}[d] & D^0 \ar@{->}[d] \\
	D^n \ar@{->}[r] & S^n \ar[r] & P^n \\
}
\]
where the square on the left is the pushout \ref{definitionsphere} and the square on the right is also a pushout. Then the outer rectangle is also a pushout. Remark that $D^n \to D^0$ is a weak equivalence and the middle vertical map is a cofibration. Since $M$ is left proper, this proves that $S^n \to P^n$ is a weak equivalence.
\epf

For a simplicial model category $M$ and $X,Y \in M$, we will write $SSet_M(X,Y)$ for the simplicial hom. Note that if $M$ is a simplicial category, then for $e \in M$, the comma category $e / M$ is also a simplicial model category. Finally, for an object of $e / M$ given by $x \in M$ and a morphism $\alpha: e \to x$, we will write just $x$ when there is no ambiguity for the map $\alpha$.

\lem\label{lemmadeloopingspheres}
Let $M$ be a simplicial model category and $n \geq 0$. For all fibrant $X \in D^0/M$, we have a weak equivalence
\[
\Omega^n SSet_{D^0/M} \left(S^0, X\right) \sim SSet_{D^0/M} \left(S^n, X\right)
\]
\endlem

\pf
For $n \geq 1$, the pushout \ref{definitionsphere} becomes a pushout in $D^0/M$ using maps such as in \ref{mapsouthpole}. Since all maps are cofibrations, this pushout is also a homotopy pushout. By applying the functor $SSet_{D^0/M}(-,X)$ to this homotopy pushout \ref{definitionsphere}, we get a homotopy pullback
\begin{equation}\label{homotopypullbackmappingspacesset}
\xymatrix{
	SSet_{D^0/M}(S^k,X) \ar[r] \ar[d] & SSet_{D^0/M}(D^k,X) \ar[d] \\
	SSet_{D^0/M}(D^k,X) \ar[r] & SSet_{D^0/M}(S^{k-1},X)
}
\end{equation}

Observe that the map $D^k \to D^0$ is a weak equivalence between cofibrant objects. Since $X$ is fibrant, as is well-known \cite[Corollary 9.3.3]{hirschhorn}, $D^k \to D^0$ induces a weak equivalence 
\[
SSet_{D^0/M}(D^0,X) \xrightarrow{\sim} SSet_{D^0/M}(D^k,X)
\]
Since $SSet_{D^0/M}(D^0,X)$ is contractible as $D^0$ is the initial object in $D^0/M$, we deduce that $SSet_{D^0/M}(D^k,X)$ is also contractible. We can therefore deduce from the homotopy pullback \ref{homotopypullbackmappingspacesset} that there is a weak equivalence
\[
\Omega SSet_{D^0/M}(S^{k-1},X) \sim SSet_{D^0/M}(S^k,X)
\]

We get the desired result by iterating this delooping.
\epf

For a simplicial model category $M$, we will write $Map_M(-,-)$ for the homotopy mapping spaces in $M$. Recall that, for $X,Y \in M$, $Map_M(X,Y)$ can be computed as $SSet_M \left(X^c,Y^f\right)$ where $X^c$ and $Y^f$ are cofibrant and fibrant replacements of $X$ and $Y$ respectively \cite{hirschhorn}.

\thm\label{generaldeloopingtheorem}
Let $M$ be a left proper simplicial model category and $n \geq 0$. For all $X \in D^0/M$, we have a weak equivalence
\[
\Omega^n Map_M \left(D^0,g_0^*(X)\right) \sim Map_{S^{n-1}/M} \left(D^n,g_n^*(X)\right)
\]
\endthm

\pf
Since $g_n^*$ preserves weak equivalences, we can assume without loss of generality that $X$ is fibrant. By adjunction and using Lemma \ref{lemmaweakequivalencespheres}, we get a weak equivalence
\[
Map_{S^{n-1}/M} \left(D^n, g_n^* (X)\right) \sim SSet_{D^0/M} \left(S^n, X\right)
\]
We get the conclusion using Lemma \ref{lemmadeloopingspheres}.
\epf

\subsection{Left cofinal Quillen functors}\label{subsectionapplicationcofinality}

\begin{definition}\label{definitionfunctorpreseringcofibrantreplacementsofterminalobject}
	A left Quillen functor $G: B \to C$ is \emph{cofinal} if $\mathbb{L}G(1)$ is contractible, where $\mathbb{L}G$ is the left derived functor of $G$ and $1$ is the terminal object in $B$.
\end{definition}


\begin{remark}\label{remarkhomotopicallycofinalandcofinalleftadjoint}
	We deduce from Lemma \ref{lemmaformulaleftderivedfunctor}  that a morphism of polynomial monads $f: S \to T$ is homotopically cofinal if and only if $f_!$ is a left cofinal Quillen functor.
\end{remark}

\begin{lemma}\label{leftquillenfunctormappingspaces}
	If $G: B \to C$ is a left cofinal Quillen functor and its right adjoint $V$ preserves weak equivalences, then for all $X \in C$ equipped with a map from the terminal object, there is a weak equivalence
	\[
		Map_C (1,X) \to Map_B (1,V(X)).
	\]
\end{lemma}

\pf
Let $1^c$ and $X^f$ be a cofibrant and fibrant replacement of the terminal object in $B$ and $X$ respectively. We have the following weak equivalences:
\[
	Map_C(1,X) \sim SSet_C(G(1^c),X^f) \sim SSet_B(1^c,V(X^f)) \sim Map_B(1,V(X)).
\]
\epf

Let $\mathrm{Bimod}_{\nfourdots}$ be the polynomial monad for quadruples $(A,B,C,D)$ where $A$ and $B$ are non-symmetric operads and $C$ and $D$ are pointed $A-B$-bimodules. Let $\mathcal{B}$ be the category of simplicial algebras of $\mathrm{Bimod}_{\nfourdots}$ and $\Phi : \mathcal{B}^{op} \to \CAT$ be the functor which sends a quadruple $(A,B,C,D)$ to the category of pointed infinitesimal $C-D$-bimodules. 
Recall that the \emph{Grothendieck construction} $\int \Phi$ is the category whose objects are pairs $(b,x)$ where $b \in \mathcal{B}$ and $x \in \Phi(b)$ and morphisms $(b,x) \to (c,y)$ are pairs $(f,\gamma)$ where $f: b \to c$ is a morphism in $\mathcal{B}$ and $\gamma: x \to \Phi(f)(y)$ is a morphism in $\Phi(b)$. There is a categorical equivalence between $\int \Phi$ and the category of simplicial algebras of the polynomial monad $\wbimodbfive$ defined in Subsection \ref{subsectiondescriptionpolynomialmonads}.

\lem
For any $b \in \mathcal{B}$, $\Phi(b)$ has a model structure.
\endlem

\pf
First note that $\Phi$ is \emph{homotopically structured} \cite{bwdgrothendieck}, that is, $\mathcal{B}$ is equipped with two classes of morphisms called \emph{horizontal} weak equivalences and fibrations, and for each $b \in \mathcal{B}$, the category $\Phi(b)$ is also equipped with classes of weak equivalences and fibrations called \emph{vertical}. Also, $\int \Phi$ is complete, cocomplete and admits a \emph{global} model structure, that is, where a morphism $(f,\gamma)$ in $\int \Phi$ is a weak equivalence (resp. fibration) if $f$ is a horizontal and $\gamma$ is a vertical weak equivalence (resp. fibration). Finally, for all $b \in \mathcal{B}$, $\Phi(b)$ is complete and cocomplete, and for all morphisms $f: b \to c$ in $\mathcal{B}$, the functor $\Phi(f)$ preserves weak equivalences, fibrations and terminal objects. It was proved in \cite{bwdgrothendieck} that, with such assumptions, $\Phi(b)$ admits a model structure for all $b \in \mathcal{B}$.
\epf

We will write $NOp$ for the category of simplicial non-symmetric operads.

\lem\label{lemmaleftcofinalfunctor}
For $b \in \mathcal{B}$ a cofibrant replacement of the terminal object, there is a left cofinal Quillen functor $\Phi(b) \to S^1 / NOp$, where $S^1$ is a circle constructed as in Definition \ref{definitionspheresanddisks}.
\endlem

\pf
Let $\mathrm{NOp}_{\bempty}$ be the polynomial monad for diagrams
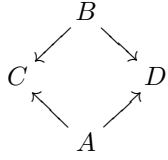
\[
\xymatrix@C=1pc@R=1pc{
	& B \ar[rd] \ar[ld] \\
	C && D \\
	& A \ar[ru] \ar[lu]
}
\]
of non-symmetric operads. There is a morphism of polynomial monads
\[
g: \mathrm{Bimod}_{\nfourdots} \to \mathrm{NOp}_{\bempty}
\]
which is given by restriction to $\{A,B,C,D\} \times \mathbb{N}$ of the morphism
\[
f: \wbimodbfive \to \nopbfive
\]
of Theorem \ref{theoremcofinality}. The classifier induced by $g$ is also the restriction to $\{A,B,C,D\} \times \mathbb{N}$ of the classifier induced by $f$. Therefore, since  $f$ is homotopically cofinal, so is $g$. Let $\mathcal{C}$ be the category of simplicial algebras of $\mathrm{NOp}_{\bempty}$ and $\Psi : \mathcal{C}^{op} \to \mathrm{CAT}$ be the functor which sends a diagram $c$ to the comma category $\colim (c) / NOp$. There is a categorical equivalence between $\int \Psi$ and the category of simplicial algebras of the polynomial monad $\nopbfive$. For $(c,y) \in \int \Psi$, the restriction functor induced by $f$ is given by
\[
f^* (c,y) = (g^*(c), W_c(y)),
\]
where $W_c: \Psi(c) \to \Phi g^*(c)$ is the obvious forgetful functor. Let $\eta$ be the unit of the adjunction induced by $g$ and, for $b \in \mathcal{B}$, let $H_b : \Phi(b) \to \Psi g_! (b)$ be the left adjoint of $\Phi(\eta_b) W_{g_!(b)}$. By adjunction, for all $(b,x) \in \int \Phi$, we have
\[
f_!(b,x) = (g_!(b),H_b(x)).
\]
From Remark \ref{remarkhomotopicallycofinalandcofinalleftadjoint}, we deduce that $f_!$ and $g_!$ are left cofinal Quillen functors. If $b \in \mathcal{B}$ is a cofibrant replacement of the terminal object, then $H_b$ is left cofinal. The diagram $g_!(b)$ forms a cofibrant replacement of the terminal object in $\mathcal{C}$, so its colimit is a circle $S^1$ in $NOp$, as constructed in Definition \ref{definitionspheresanddisks}. This means that the categories $\Psi g_!(b)$ and $S^1 / NOp$ are equivalent, and $H_b$ gives us the desired functor.
\epf

\subsection{Fibration sequence theorem}\label{subsectionfibrationsequence}

\lem\label{lemmahomotopypushoutcoproduct}
Let $M$ be a model category. Let $x \to a \to b$ be a composite in $M$. Assume $x$ cofibrant and the map $a \to b$ is a cofibration. The following square is a homotopy pushout in $M$
\begin{equation}\label{homotopypushout}
\xymatrix{
	a \amalg x \ar[r] \ar[d] & a \ar[d] \\
	b \amalg x \ar[r] & b
}
\end{equation}
\endlem

\pf
Let $a^c$ be an object in $M$ together with a factorisation
\[
\xymatrix{
	a \amalg x \ar[rr] \ar@{ >->}[rd] && a \\
	& a^c \ar[ru]_\sim
}
\]
Let us consider the following commutative diagram
\[
\xymatrix{
	a \ar@{ >->}[r] \ar@{ >->}[d] \ar@{}[dr] | {i} & a \amalg x \ar@{ >->}[r] \ar@{ >->}[d] \ar@{}[dr] | {ii} & a^c \ar@{ >->}[d] \\
	b \ar@{ >->}[r] & b \amalg x \ar@{ >->}[r] & p
}
\]
where $p$ is defined as the pushout of the square $ii$.

We want to prove that the composite $b \to p$ is a weak equivalence. It is actually a trivial cofibration. To see this, observe first that the composite $a \to a^c$ is a trivial cofibration. Indeed, it is a weak equivalence by the two-out-of-three property and it is also a cofibration since $x$ is cofibrant. Moreover, since the square $i$ is a pushout, the big rectangle is also a pushout. We get the conclusion using the fact that trivial cofibrations are stable under pushout.

Since the square $ii$ is a pushout and all the maps in this square are cofibrations, it is also a homotopy pushout. Since all the objects in the square $ii$ are weakly equivalent to the objects in the square \ref{homotopypushout}, this last square is also a homotopy pushout.
\epf

\thm\label{lemmaforfibrationsequence}
Let $e \in M$ be cofibrant. Then for all maps $r \to s$ in $e/M$, there is a fibration sequence
\[
Map_{e/M} (r,s) \to Map_M (V(r),V(s)) \to Map_M (e,V(s)),
\]
where $V: e/M \to M$ is the forgetful functor. Moreover, if $M$ is left proper, the assumption that $e$ is cofibrant can be dropped.
\endthm

\pf

Since $V$ preserves weak equivalences, we can assume without loss of generality that $r$ is cofibrant and $s$ is fibrant in $e/M$. By definition of being cofibrant in $e/M$, there is a cofibration $e \to V(r)$. Since $e$ is cofibrant, this means that $V(r)$ is also cofibrant in $M$. Let $G$ be the left adjoint of $V$. Then
\[
VG(x) = x \amalg e.
\]
Let $\iota$ be the initial object in $e/M$, given by the identity map on $e$. We have the following commutative square in $e/M$:
\begin{equation}\label{homotopypushoutforfibrationsequence}
\xymatrix{
	G(e) \ar[r] \ar[d] & \iota \ar[d] \\
	GV(r) \ar[r] & r
}
\end{equation}
Observe that a square in $e/M$ is a homotopy pushout if and only if $V$ applied to this square gives a homotopy pushout. But when we apply $V$ to the square \ref{homotopypushoutforfibrationsequence}, we get the square
\[
\xymatrix{
	e \amalg e \ar[r] \ar[d] & e \ar[d] \\
	V(r) \amalg e \ar[r] & V(r)
}
\]
which is a homotopy pushout by Lemma \ref{lemmahomotopypushoutcoproduct}. This proves that the square \ref{homotopypushoutforfibrationsequence} is a homotopy pushout. Applying $SSet_{e/M}(-,s)$ to \ref{homotopypushoutforfibrationsequence} gives us the homotopy pullback
\[
\xymatrix{
	SSet_{e/M}(r,s) \ar[r] \ar[d] & SSet_{e/M}(\iota,s) \ar[d] \\
	SSet_{e/M}(GV(r),s) \ar[r] & SSet_{e/M}(G(e),s)
}
\]
By adjunction, the bottom map is equivalent to the map
\[
SSet_M(V(r),V(s)) \to SSet_M(e,V(s))
\]
Since $\iota$ is the initial object in $e/M$, $SSet_{e/M}(\iota,s)$ is contractible.

If $M$ is left proper, any weak equivalence $e \to e'$ induces a Quillen equivalence between $e / M$ and $e' / M$ \cite[Proposition 2.7]{rezk}, so we can drop the assumption that $e$ is cofibrant. This concludes the proof.
\epf

For a functor $U: B \to M$ and $e \in M$, we write $e/U$ for the \emph{slice construction}, whose objects are pairs $(b,\alpha)$, with $b \in B$ and $\alpha: e \to U(b)$ a morphism in $M$.

\thm\label{theoremfibrationsequence}
Let $U : B \to M$ be a right Quillen functor and $e \in M$ cofibrant. Then for all maps $r \to s$ in $e/U$, there is a fibration sequence
\[
Map_{e/U} (r,s) \to Map_B (V(r),V(s)) \to Map_M (e,UV(s))
\]
where $V: e/U \to B$ is the forgetful functor.
\endthm

\pf
Let $F: M \to B$ be the left adjoint of $U$. Since $F(e)$ is cofibrant in $B$, we can apply Lemma \ref{lemmaforfibrationsequence} to get, for all maps $r \to s$ in $F(e) / B$, the fibration sequence
\[
Map_{F(e) / B} (r,s) \to Map_B (V(r),V(s)) \to Map_B (F(e),V(s))
\]

We get the conclusion by adjunction and from the fact that the categories $F(e) / B$ and $e / U$ are equivalent.
\epf

\subsection{Equivalence between nerves of classifiers}

Recall from Subsection \ref{subsectionapplicationcofinality} that $\mathrm{Bimod}_{\nfourdots}$ is the polynomial monad for quadruples $(A,B,C,D)$ where $A$ and $B$ are non-symmetric operads and $C$ and $D$ are pointed $A-B$-bimodules. We write $\mathcal{B}$ for the category of algebras of $\mathrm{Bimod}_{\nfourdots}$. Also recall from Subsection \ref{subsectiondescriptionpolynomialmonads} that $\wbimodbfive$ is the polynomial monad for quintuples $(A,B,C,D,E)$ where $(A,B,C,D)$ is an element of $\mathcal{B}$ and $E$ is a pointed infinitesimal $C-D$-bimodule. There is a morphism of polynomial monads
\begin{equation}\label{classifier1}
\mathrm{Bimod}_{\nfourdots} \to \wbimodbfive
\end{equation}
such that the restriction functor forgets the pointed infinitesimal $C-D$-bimodule $E$.

Now let $\mathrm{BimodSP}_{\nfivedots}$ be the polynomial monad for quintuples $(A,B,C,D,E)$ where $(A,B,C,D)$ is an element of $\mathcal{B}$ and $E$ is a pair $(E_0,E_1)$ of simplicial sets equipped with a pair of maps
\begin{equation}\label{pairofmaps}
	A_2 \otimes E_0 \rightrightarrows E_1.
\end{equation}
Let $\wbimodbfive^\circ$ be the polynomial monad for quintuples $(A,B,C,D,E)$ where $(A,B,C,D)$ is an element of $\mathcal{B}$ and $E$ is an infinitesimal $C-D$-bimodule. 
There is a morphism of polynomial monads
\begin{equation}\label{classifier2}
\mathrm{BimodSP}_{\nfivedots} \to \wbimodbfive^\circ
\end{equation}
such that the restriction functor forgets the infinitesimal $C-D$-bimodule $E$ except for the pair $(E_0,E_1)$ and the maps \ref{pairofmaps} are given by the composites
\begin{equation}\label{composite1}
A_2 \otimes E_0 \simeq A_2 \otimes E_0 \otimes I \xrightarrow{1 \otimes 1 \otimes \zeta} A_2 \otimes E_0 \otimes D_1 \xrightarrow{\phi} E_1
\end{equation}
and
\begin{equation}\label{composite2}
A_2 \otimes E_0 \simeq A_2 \otimes I \otimes E_0 \xrightarrow{1 \otimes \zeta \otimes 1} A_2 \otimes C_1 \otimes E_0 \xrightarrow{\phi} E_1,
\end{equation}
where $\zeta$ is the map given by the square \ref{axiompointedbimodule} and $\phi$ is the map \ref{mapinfinitesimalleft}. The objective of this subsection is to prove the following lemma:
\begin{lemma}\label{lemmaweakequivalenceclassifiers}
	There is a weak equivalence between the nerve of the classifier induced by the morphism \ref{classifier1} and the nerve of the classifier induced by the morphism \ref{classifier2}.
\end{lemma}

\pf
We define two new morphisms of polynomial monads very similar to the ones defined above. Let $\mathrm{Bimod}_{\nspecialthreedots}$ be the polynomial monad whose algebras are the same as the algebras of $\mathrm{Bimod}_{\nfourdots}$, but we force the non-symmetric operad $A$ to be the terminal non-symmetric operad $\mathcal{A}ss$. That is, the polynomial monad for quadruples $(A,B,C,D)$ where $A=\mathcal{A}ss$, $B$ is any non-symmetric operad and $C$ and $D$ are pointed $A-B$-bimodules. In a totally similar way, one can define the polynomial monads $\mathrm{IBimod}_{\nspecialfourdots}$, $\mathrm{BimodSP}_{\nspecialfourdots}$ and $\mathrm{IBimod}_{\nspecialfourdots}^\circ$. We have two more morphisms of polynomial monads
\begin{equation}\label{classifier3}
\mathrm{Bimod}_{\nspecialthreedots} \to \mathrm{IBimod}_{\nspecialfourdots}
\end{equation}
and
\begin{equation}\label{classifier4}
\mathrm{BimodSP}_{\nspecialfourdots} \to \mathrm{IBimod}_{\nspecialfourdots}^\circ
\end{equation}

To simplify the notations, let us fix an integer $n \geq 0$ and let us write $\mathcal{X}$, $\mathcal{X}'$, $\mathcal{Y}$ and $\mathcal{Y}'$ for the underlying categories indexed by the label $E$ and the integer $n$ of the classifiers induced by the morphisms \ref{classifier2}, \ref{classifier4}, \ref{classifier1} and \ref{classifier3} respectively. Let us describe these categories more explicitly. For the sake of brevity, we will only describe the objects and generating morphisms of each category. We leave it to the reader to figure out the obvious relations between the generating morphisms and to check that the description of each category corresponds to the right classifier. The objects of $\mathcal{Y}$ are isomorphism classes of planar trees whose vertices have labels in $\{A,B,C,D\}$ lying on a line and a point as defined in \ref{conditionstar}. So there can be no vertices labelled with $E$. The generating morphisms are contractions of subtrees with vertices only labelled with $A$ (resp. $B$) to a vertex labelled with $A$ (resp. $B$) and contraction of subtrees which do not contain vertices with label $D$ (resp. $C$) to a vertex with label $C$ (resp. $D$). The objects of $\mathcal{X}$ are isomorphism classes of planar trees whose vertices have labels in $\{A,B,C,D,E\}$ lying on a line and a point. There must be exactly one vertex with label $E$ and this vertex must have zero or one edge above it. The generating morphisms are the same as $\mathcal{Y}$, plus contraction of the following subtrees to an unary vertex with label $E$, assuming that there are no vertices labelled with $C$ or $D$ above it:
\[
	\begin{tikzpicture}[scale=.5]
		\draw (0,-.5) -- (0,0) -- (-1,1);
		\draw (0,0) -- (.7,.7);
		\draw[fill] (0,0) circle (1.5pt) node[left]{$A$};
		\draw[fill] (.7,.7) circle (1.5pt) node[above]{$E$};
		\draw (6,-.5) -- (6,0) -- (7,1);
		\draw (6,0) -- (5.3,.7);
		\draw[fill] (6,0) circle (1.5pt) node[left]{$A$};
		\draw[fill] (5.3,.7) circle (1.5pt) node[above]{$E$};
	\end{tikzpicture}
\]
The objects of $\mathcal{Y}'$ are linearly ordered finite sets, possibly empty, of isomorphism classes of planar trees, which we will call \emph{forests}. Each tree in the forest has vertices labelled in $\{B,C,D\}$. Only the root vertices may have a label which is not $B$, all the other vertices must have label $B$. All the trees whose root vertex is labelled with $C$ must be to the left of all trees whose root vertex is labelled with $D$. The generating morphisms are the same as those for $\mathcal{Y}$ which do not involve vertices with label $A$. Another type of generating morphism is the concatenation of trees next to each other. It is given by gluing the root vertices of each tree together to get a unique root vertex as in the following picture:
\[
	\begin{tikzpicture}[scale=.6]
		\draw (0,-.5) -- (0,0) -- (-1,1);
		\draw (0,1) -- (0,0) -- (1,1);
		\draw (3,-.5) -- (3,0) -- (2.3,1);
		\draw (3,0) -- (3.7,1) -- (3.2,2);
		\draw (3.7,1) -- (4.2,2);
		\draw[fill] (0,0) circle (1.2pt) node[left]{$C$};
		\draw[fill] (0,1) circle (1.2pt) node[above]{$B$};
		\draw[fill] (3,0) circle (1.2pt) node[left]{$C$};
		\draw[fill] (3.7,1) circle (1.2pt) node[right]{$B$};
		\draw (7,0) node[above]{$\longrightarrow$};
		\draw (11,-.5) -- (11,0) -- (9.6,1);
		\draw (10.3,1) -- (11,0) -- (11,1);
		\draw (11.7,1) -- (11,0) -- (12.4,1) -- (11.9,2);
		\draw (12.4,1) -- (12.9,2);
		\draw[fill] (11,0) circle (1.2pt) node[left]{$C$};
		\draw[fill] (10.3,1) circle (1.2pt) node[above]{$B$};
		\draw[fill] (12.4,1) circle (1.2pt) node[right]{$B$};
	\end{tikzpicture}
\]
The trees next to each other must all have their root vertices labelled with $C$ (resp. $D$) and the obtained concatenated tree will have its root vertex labelled with $C$ (resp. $D$). Finally, regarding $\mathcal{X}'$, its objects are forests where each vertex has a label in $\{B,C,D,E\}$. The conditions are as for $\mathcal{Y}'$ but there must also be a unique vertex with label $E$, which satisfies the condition of lying on a line and a point. So all the vertices labelled with $C$ must be to the left of the vertex labelled with $E$, and all the vertices labelled with $D$ to the right of it. The unique vertex with label $E$ must be either the vertex of a trunk, that is a corolla without leaves, or the unary root vertex of a tree whose other vertices all have label $B$. The generating morphisms are as described for $\mathcal{Y}'$. There are also the generating morphisms which remove the trunk whose vertex is labelled with $E$ and adds a unary vertex with label $E$ to the root edge of the tree situated to the left or to the right of where the trunk was originally.

Let us now construct a diagram of functors
\begin{equation}\label{diagramsmoothfunctors}
\xymatrix{
	\mathcal{X} \ar[r]^-F & \mathcal{X}' \ar[r]^-G & \mathcal{Y}' & \mathcal{Y}. \ar[l]_-H
}
\end{equation}
The functors $F$ and $H$ cuts all the edges above vertices labelled with $A$ and keeps only the branches lying just above these vertices. The functor $G: \mathcal{X}' \to \mathcal{Y}'$ acts on a tree as follows. If the unique vertex with label $E$ has no incoming edges, then it removes the trunk. If the unique vertex with label $E$ has an incoming edge, then it removes this vertex.
%
%
%
%
%
%
%
%
Such functors on objects can obviously be extended on morphisms.

First note that the fibres of each functor have contractible nerve. Indeed, the fibres of $F$ and $H$ have a terminal object. 
%
The fibres of $G$ consist of zigzags of morphisms. For example, the fibre over the forest	
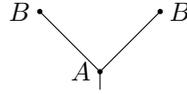
\[
\begin{tikzpicture}[scale=.8]
\draw (0,0) -- (0,1);
\draw (.7,0) -- (.7,1);
\draw[fill] (0,1) circle (1pt) node[above]{$\tiny{B}$};
\draw[fill] (.7,1) circle (1pt) node[above]{$\tiny{B}$};
%
\end{tikzpicture}
\]
is given by:
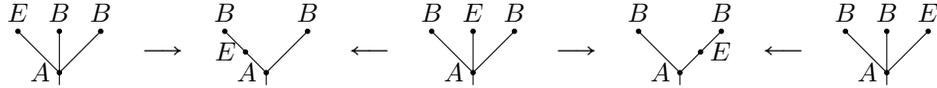
\[
\begin{tikzpicture}[scale = .8]
\draw (-.7,0) -- (-.7,1);
\draw (0,0) -- (0,1);
\draw (.7,0) -- (.7,1);

\draw[fill] (-.7,1) circle (1pt) node[above]{$\tiny{E}$};
\draw[fill] (0,1) circle (1pt) node[above]{$\tiny{B}$};
\draw[fill] (.7,1) circle (1pt) node[above]{$\tiny{B}$};

\draw (2,.5) node[]{$\longrightarrow$};

\begin{scope}[shift={(4,0)}]
\draw (-.35,0) -- (-.35,1);
\draw (.35,0) -- (.35,1);

\draw[fill] (-.35,.5) circle (1pt) node[left]{$\tiny{E}$};
\draw[fill] (-.35,1) circle (1pt) node[above]{$\tiny{B}$};
\draw[fill] (.35,1) circle (1pt) node[above]{$\tiny{B}$};
\end{scope}

\draw (6,.5) node[]{$\longleftarrow$};

\begin{scope}[shift={(8,0)}]
\draw (-.7,0) -- (-.7,1);
\draw (0,0) -- (0,1);
\draw (.7,0) -- (.7,1);

\draw[fill] (-.7,1) circle (1pt) node[above]{$\tiny{B}$};
\draw[fill] (0,1) circle (1pt) node[above]{$\tiny{E}$};
\draw[fill] (.7,1) circle (1pt) node[above]{$\tiny{B}$};
\end{scope}

\draw (10,.5) node[]{$\longrightarrow$};

\begin{scope}[shift={(12,0)}]
\draw (-.35,0) -- (-.35,1);
\draw (.35,0) -- (.35,1);

\draw[fill] (.35,.5) circle (1pt) node[right]{$\tiny{E}$};
\draw[fill] (-.35,1) circle (1pt) node[above]{$\tiny{B}$};
\draw[fill] (.35,1) circle (1pt) node[above]{$\tiny{B}$};
\end{scope}

\draw (14,.5) node[]{$\longleftarrow$};

\begin{scope}[shift={(16,0)}]
\draw (-.7,0) -- (-.7,1);
\draw (0,0) -- (0,1);
\draw (.7,0) -- (.7,1);

\draw[fill] (-.7,1) circle (1pt) node[above]{$\tiny{B}$};
\draw[fill] (0,1) circle (1pt) node[above]{$\tiny{B}$};
\draw[fill] (.7,1) circle (1pt) node[above]{$\tiny{E}$};
\end{scope}

\end{tikzpicture}
\]
Let us now prove that each functor is smooth, which will conclude the proof thanks to Cisinski lemma \ref{cisinskilemma}. For the functor $F$, let us take an object $x_1$ in $\mathcal{X}$ and a map $f_1: y_0 \to y_1$ in $\mathcal{X}'$ such that $F(x_1)=y_1$. The lifting category of $f_1$ over $x_1$ has a terminal object, which can be constructed in the following way. From the tree $x_1$, one can get a tree $T$, whose vertices are all labelled with $A$, and which is obtained by making the same cuts as for the definition of the functor $F$, but keeping the tree below instead of the forest above. There are trees $T_1,\ldots,T_k$, with $k$ the number of leaves of $T$, such that the tree obtained by attaching these trees to each leaf of $T$ gives an object in the fibre of $y_0$. Moreover, one can take $T_1,\ldots,T_k$ with at most one vertex labelled with $A$, which is not unary. Then the obtained tree is the terminal object of the lifting category. Similarly, $H$ is also a smooth functor. For the functor $G$, let us take an object $x_1$ in $\mathcal{X}'$ and a map $f_1: y_0 \to y_1$ in $\mathcal{Y}'$ such that $G(x_1) = y_1$.
The lifting category of a map $f_1$ over $x_1$ is the terminal category if the unique vertex labelled with $E$ in $x_1$ has $0$ incoming edges and consists in a cospan if it has $1$ incoming edge.
\epf

%


\subsection{Direct double delooping proof}

Let $\mathcal{B}$ be the category and $\Phi: \mathcal{B}^{op} \to \mathrm{CAT}$ be the functor defined in the previous subsection. Let $\Phi^\circ: \mathcal{B}^{op} \to \mathrm{CAT}$ be the functor which sends a quadruple $(A,B,C,D)$ to the category of infinitesimal $C-D$-bimodules. Finally, let $\mathfrak{b}$ be a cofibrant replacement of the terminal object in $\mathcal{B}$ and $V: \Phi(\mathfrak{b}) \to \Phi^\circ(\mathfrak{b})$ be the forgetful functor.

\lem\label{tdh3}
Let $x \to y$ be a map in $\Phi(\mathfrak{b})$. If $y_0$ and $y_1$ are contractible, then there is a weak equivalence
\begin{equation}\label{we3}
Map_{\Phi(\mathfrak{b})} (x,y) \to Map_{\Phi^\circ(\mathfrak{b})} (V(x),V(y)).
\end{equation}
\endlem

\pf
Let $\Psi: \mathcal{B}^{op} \to \mathrm{CAT}$ be the functor which sends a quadruple $(A,B,C,D)$ to the category of pairs $(E_0,E_1)$ of simplicial sets equipped with a pair of maps \ref{pairofmaps}. Let $U: \Phi^\circ(\mathfrak{b}) \to \Psi(\mathfrak{b})$ be the forgetful functor which sends an infinitesimal bimodule $E$ to the pair $(E_0,E_1)$ together with the maps \ref{pairofmaps} given by the composites \ref{composite1} and \ref{composite2}. Finally, let $e$ be a cofibrant replacement of the terminal object in $\Psi(\mathfrak{b})$ and $x \to y$ be a map in $e/U$. Applying Theorem \ref{theoremfibrationsequence}, we get a fibration sequence
\[
Map_{e/U} (x,y) \to Map_{\Phi^\circ(\mathfrak{b})} (V(x),V(y)) \to Map_{\Psi(\mathfrak{b})} (e,UV(y)).
\]
If $y_0$ and $y_1$ are contractible, then the last mapping space in this fibration sequence is also contractible, which means that the first map is a weak equivalence. It remains to show that there is a Quillen equivalence between $e/U$ and $\Phi(\mathfrak{b})$.

Let $F$ be the left adjoint of $U$. Then the category $e/U$ is isomorphic to $F(e)/\Phi^\circ(\mathfrak{b})$. We now prove that there is a weak equivalence between $F(e)$ and $V(0)$, where $0$ is the initial object of $\Phi(\mathfrak{b})$. Let us denote by $h$ and $k$ the maps \ref{classifier1} and \ref{classifier2} respectively. According to Lemma \ref{lemmaformulaleftderivedfunctor}, the classifiers induced by these maps are $h_!(\mathfrak{b})$ and $k_!(\mathfrak{b},e)$ respectively. On the other hand, $h_!(\mathfrak{b}) = (\mathfrak{b},0)$ and $k_!(\mathfrak{b},e) = (\mathfrak{b},F(e))$. The desired weak equivalence can therefore be deduced from Lemma \ref{lemmaweakequivalenceclassifiers}. The category $\Phi^\circ(\mathfrak{b})$, as a category of presheaves, is left proper. This implies \cite[Proposition 2.7]{rezk} that the weak equivalence between $F(e)$ and $V(0)$ induces a Quillen equivalence between $F(e)/\Phi^\circ(\mathfrak{b})$ and $V(0)/\Phi^\circ(\mathfrak{b})$. Finally, note that equipping an infinitesimal bimodule $E$ with a map $V(0) \to E$ corresponds to making it pointed, in other words, the categories $V(0)/\Phi^\circ(\mathfrak{b})$ and $\Phi(\mathfrak{b})$ are isomorphic. This concludes the proof.
\epf

\pfof[Theorem \ref{theoremturchindwyerhess}]
Since the category of simplicial non-symmetric operads $NOp$ is left proper \cite{bataninberger,muro}, we can apply Theorem \ref{generaldeloopingtheorem} to get a weak equivalence
\begin{equation}\label{we1}
\Omega^2  Map_{NOp} (D^0, g_0^* (y)) \to Map_{S^1/NOp} (D^2, g_2^* (y)),
\end{equation}
where disks and spheres are constructed as in Definition \ref{definitionspheresanddisks} and $g_n: S^{n-1} \to D^0$ is the map \ref{mapgn}. Let $H_\mathfrak{b} : \Phi(\mathfrak{b}) \to S^1 / NOp$ be the left cofinal Quillen functor given by Lemma \ref{lemmaleftcofinalfunctor}, and $W_\mathfrak{b}$ its right adjoint. By adjunction, for $x \in \Phi(\mathfrak{b})$ cofibrant and $y \in S^1 / NOp$ fibrant, there is an isomorphism
\begin{equation}\label{we2}
Map_{S^1/NOp} (H_\mathfrak{b}(x),y) \to Map_{\Phi(\mathfrak{b})} (x,W_\mathfrak{b}(y)).
\end{equation}
The assumption that $y$ is fibrant can be dropped since $W_\mathfrak{b}$ preserves weak equivalences. Also, if $x$ is a cofibrant replacement of the terminal object, $H_\mathfrak{b}(x)$ is contractible since $H_\mathfrak{b}$ is left cofinal. Let $1$ be the terminal object in $\mathcal{B}$ and again in $\Phi^\circ(1)$, abusing notations. The unique map $\mathfrak{b} \to 1$ is a weak equivalence, so it induces a Quillen equivalence between the categories of simplicial presheaves $\Phi^\circ(\mathfrak{b})$ and $\Phi^\circ(1)$. Therefore, putting together the weak equivalences \ref{we1}, \ref{we2} and \ref{we3}, we get, for a multiplicative operad $\mathcal{O}$ such that $\mathcal{O}_0$ and $\mathcal{O}_1$ are contractible, a weak equivalence
\[
\Omega^2  Map_{NOp} (\mathcal{A}ss,u^*(\mathcal{O})) \sim Map_{\Phi^\circ(1)} (1,f^*(\mathcal{O})).
\]
According to Lemma \ref{theoremequivalenceweakbimodules}, $\Phi^\circ(1)$ is equivalent to the category of infinitesimal $Ass$-bimodules in the classical sense, which is in turn equivalent to the category cosimplicial objects \cite[Lemma 4.2]{turchincosimplicial}. This concludes the proof.
\epf

\noindent {\bf Acknowledgements.} I wish to express my  gratitude to Michael Batanin and Steve Lack for many useful discussions, and Victor Turchin, Marcy Robertson and Boris Shoikhet for reading, commenting and giving me feedback on my PhD thesis. I also would like to thank the anonymous referee for their remarks and suggestions of improvement.

\bibliographystyle{plain}
\bibliography{references}

\begin{thebibliography}{10}

\bibitem{baezdolan}
John~C Baez and James Dolan.
\newblock Higher-dimensional algebra iii. n-categories and the algebra of
  opetopes.
\newblock {\em Advances in Mathematics}, 135(2):145--206, 1998.

\bibitem{batanin}
Michael Batanin.
\newblock The {E}ckmann--{H}ilton argument and higher operads.
\newblock {\em Advances in Mathematics}, 217(1):334--385, 2008.

\bibitem{bataninberger}
Michael Batanin and Clemens Berger.
\newblock Homotopy theory for algebras over polynomial monads.
\newblock {\em Theory and Applications of Categories}, 32(6):148--253, 2017.

\bibitem{batanindeleger}
Michael Batanin and Florian De~Leger.
\newblock Polynomial monads and delooping of mapping spaces.
\newblock {\em Journal of Noncommutative Geometry}, 13(4):1521+, 2019.

\bibitem{bwdgrothendieck}
Michael Batanin, Florian De~Leger, and David White.
\newblock Model structures on operads and algebras from a global perspective.
\newblock {\em arXiv:2311.07320}, 2023.

\bibitem{boavida}
Pedro Boavida~de Brito and Michael Weiss.
\newblock Spaces of smooth embeddings and configuration categories.
\newblock {\em Journal of Topology}, 11(1):65--143, 2018.

\bibitem{cisinski}
Denis-Charles Cisinski.
\newblock {\em Les pr{\'e}faisceaux comme mod{\`e}les des types d'homotopie}.
\newblock Soci{\'e}t{\'e} math{\'e}matique de France, 2006.

\bibitem{delegergrego}
Florian De~Leger and Maro{\v{s}} Grego.
\newblock Triple delooping for multiplicative hyperoperads.
\newblock {\em arXiv:2309.15055}, 2023.

\bibitem{ducoulombier}
Julien Ducoulombier.
\newblock Delooping derived mapping spaces of bimodules over an operad.
\newblock {\em Journal of Homotopy and Related Structures}, 14(2):411--453,
  2019.

\bibitem{ducoulombierturchin}
Julien Ducoulombier and Victor Turchin.
\newblock Delooping the functor calculus tower.
\newblock {\em arXiv preprint arXiv:1708.02203}, 2017.

\bibitem{dwyerhess}
William Dwyer and Kathryn Hess.
\newblock Long knots and maps between operads.
\newblock {\em Geometry \& Topology}, 16(2):919--955, 2012.

\bibitem{gambinokock}
Nicola Gambino and Joachim Kock.
\newblock Polynomial functors and polynomial monads.
\newblock {\em Mathematical Proceedings of the Cambridge Philosophical
  Society}, 154(1):153--192, 2013.

\bibitem{hirschhorn}
Philip Hirschhorn.
\newblock {\em Model categories and their localizations}.
\newblock Number~99. American Mathematical Society, 2009.

\bibitem{markl2008operads}
Martin Markl.
\newblock Operads and props.
\newblock {\em Handbook of algebra}, 5:87--140, 2008.

\bibitem{muro}
Fernando Muro.
\newblock Homotopy theory of non-symmetric operads, ii: Change of base category
  and left properness.
\newblock {\em Algebraic \& Geometric Topology}, 14(1):229--281, 2014.

\bibitem{rezk}
Charles Rezk.
\newblock Every homotopy theory of simplicial algebras admits a proper model.
\newblock {\em Topology and its Applications}, 119(1):65--94, 2002.

\bibitem{sinha1}
Dev Sinha.
\newblock Operads and knot spaces.
\newblock {\em Journal of the American Mathematical Society}, 19(2):461--486,
  2006.

\bibitem{turchincosimplicial}
Victor Turchin.
\newblock Hodge-type decomposition in the homology of long knots.
\newblock {\em Journal of Topology}, 3(3):487--534, 2010.

\bibitem{turchin}
Victor Turchin.
\newblock Delooping totalization of a multiplicative operad.
\newblock {\em Journal of Homotopy and Related Structures}, 9(2):349--418,
  2014.

\bibitem{weber}
Mark Weber.
\newblock Internal algebra classifiers as codescent objects of crossed internal
  categories.
\newblock {\em Theory and Applications of Categories}, 30(50):1713--1792, 2015.

\end{thebibliography}

\end{document}